\documentclass[12pt,draft]{article}

\usepackage{amsthm,amssymb,amsmath,amscd,amstext,mathrsfs}
\usepackage{latexsym, mathptmx}

\newtheorem{theorem}{Theorem}[section]
\newtheorem{definition}{Definition}[section]
\newtheorem{lemma}{Lemma}[section]

\newtheorem{remark}{Remark}[section]

\newcommand{\nc}{\newcommand}

\nc{\OO}{{\mathbb O}}
\nc{\HH}{{\mathbb H}}
\nc{\C}{{\mathbb C}}
\nc{\R}{{\mathbb R}}
\nc{\Z}{{\mathbb Z}}
\nc{\N}{{\mathbb N}}

\nc{\ii}{{\bf i}}
\nc{\dd}{{\rm d}}

\nc{\cC}{{\mathscr C}}
\nc{\ce}{{\mathscr E}}
\nc{\cf}{{\mathscr F}}
\nc{\cg}{{\mathscr G}}
\nc{\ch}{{\mathscr H}}
\nc{\cm}{{\mathscr M}}
\nc{\co}{{\mathscr O}}

\begin{document}

\title{Complex structure on the six dimensional sphere\\ 
from a spontaneous symmetry breaking}

\author{G\'abor Etesi\\ 
\small{{\it Department of Geometry, Mathematical Institute, Faculty of 
Science,}}\\ 
\small{{\it Budapest University of Technology and Economics,}}\\
\small{{\it Egry J. u. 1, H \'ep., H-1111 Budapest, Hungary}}
\footnote{e-mail: {\tt etesi@math.bme.hu}}}

\maketitle

\pagestyle{myheadings}
\markright{G. Etesi: Complex structure on the six-sphere from a spontaneous 
symmetry breaking}

\thispagestyle{empty}

\centerline{\it\large Dedicated to the memory of my father}

\begin{abstract}
Existence of a complex structure on the $6$ dimensional sphere is 
proved in this paper. The proof is based on re-interpreting a hypothetical 
complex structure as a classical ground state of a Yang--Mills--Higgs-like 
theory on $S^6$. This classical vacuum solution is then constructed by 
Fourier expansion (dimensional reduction) from the obvious one of a similar 
theory on the $14$ dimensional exceptional compact Lie group ${\rm G}_2$. 
\end{abstract}

\centerline{AMS Classification: Primary: 32J17, 32J25; Secondary: 43A30, 
53C56, 81T13}

\centerline{Keywords: {\it Six-sphere; Complex structure; 
Fourier expansion; Yang--Mills--Higgs theory}}


\section{Introduction}
\label{one}


A classical result of Borel and Serre from 1951 states that among 
spheres the two and the six dimensional are the only ones which can 
carry almost complex structures \cite{bor-ser1, bor-ser2}. In case of the 
two-sphere the almost complex structure $I$ is {\it unique} and 
stems from the embedding $S^2\subset {\rm Im}\HH$ as the unit sphere hence 
$I$ acts at $x\in S^2$ simply by multiplication with $x$ itself. Equivalently 
$I$ can be constructed through the identification $S^2\cong\C P^1$ 
consequently it is {\it integrable} in the sense that it comes from a 
complex manifold structure. On the contrary the six-sphere admits a 
{\it plethora} of almost complex structures which are 
typically {\it non-integrable}. For example the analogouos six 
dimensional almost complex structure $I$ can be constructed \cite[pp. 
163-164]{wil} from the inclusion $S^6\subset {\rm Im}\OO$ where $\OO$ 
refers to the octonions or Cayley numbers; however it is not 
integrable \cite{fro}. This can be directly proved by a lengthy 
calculation of its Nijenhuis tensor: since it is not zero this Cayley almost 
complex structure is not integrable by the Newlander--Nirenberg 
theorem \cite{new-nir}. Putting an orientation and a Riemannian metric $g$ 
onto $S^6$ an abundance of other almost complex structures compatible with the 
orientation and the metric emerges as smooth sections of the projectivized 
positive chirality spinor bundle $P\Sigma^+$, or equivalently, of the positive 
twistor space $Z^+(S^6,[g])$ where $[g]$ denotes the conformal class of $g$, 
cf. \cite[Chapter IV, \S 9]{law-mic}. Homotopy theory also can be used to 
construct almost complex structures \cite{pen-tan}. In spite of the two 
dimensional case however, it has been unknown whether or not any of these or 
other almost complex structures are integrable \cite[pp. 16-18]{che1}---or 
it has rather been believed that all of them are not-integrable hence 
the six-sphere cannot be a complex manifold at all \cite[p. 424]{gre-sch-wit}.

One source of this dismissive attitude might 
be an \ae sthetic aversion: if the most canonical and natural 
Cayley almost complex structure is not integrable then what else could be? 
Another one might be the feeling that the existence of a complex structure 
on $S^6$ would be somehow awkward: for example \cite{huc-keb-pet} 
it would allow one to perform non-K\"ahlerian deformations of the standard 
complex structure on $\C P^3$; a ``minor disaster'' as LeBrun says 
in \cite{leb}. 

During the past six decades several works appeared 
which claimed to prove the non-existence however one-by-one they turned 
out to be erroneous. Examples are \cite{adl, hsi} while the latest 
trial was Chern's attack \cite{che2} in 2004 based on the link between 
$S^6$ and the exceptional group ${\rm G}_2$. Although it 
also seems to have a gap (cf. \cite{bry2}), Chern's approach has 
introduced interesting new techniques into gauge theories \cite{mar-nie}.

These failures are not surprising because this complex structure, if 
exists, cannot be grasped by conventional means. Indeed, let $X$ 
denote this hypothetical compact complex $3$-manifold i.e., an 
integrable almost complex structure $J$ on $S^6$. Regarding
classical geometry, $X$ cannot have anything to do with $S^6\subset\R^7$
as the usual sphere with the standard round metric $g_0$ since $J$ cannot 
be orthogonal with respect to it \cite{leb, zho1} or even to a metric $g$ in 
an open neighbourhood of $g_0$ \cite{zho3}. That is, $X$ cannot 
be constructed from a section of $Z^+(S^6, [g])$ in the sense above with a 
metric $g$ ``close'' to $g_0$ (for some properties of the 
standard twistor space $Z^+(S^6, [g_0])$ cf. \cite{mus, slu, zho1} and see 
also \cite{zho2}). Since $\pi_k(X)\cong 0$ ($k=1,\dots,5$) 
$X$ cannot be fibered over any lower dimensional manifold therefore 
algebraic topology does not help to cut down the problem to a 
simpler one. It also follows that $H^2(X;\Z )\cong 0$ consequently $X$ 
cannot be projective algebraic or even K\"ahler hence cannot be approached 
by conventional algebraic geometry or analysis. In fact any meromorphic 
function on $X$ must be constant \cite{cam-dem-pet} demonstrating that the 
algebraic dimension of $X$ is zero consequently it is as far from being 
algebraic as possible. Lacking good functions, the inadequacy of powerful 
methods of complex analysis also follows. Last but not least we also know that 
$X$ is a truely inhomogeneous complex manifold in the sense that 
${\rm Aut}_{\co}X$ cannot act transitively and none of its orbits can be 
open \cite{huc-keb-pet} (actually this property would permit the existence of 
the aforementioned exotic complex projective spaces). 

Nevertheless the Hodge numbers of $X$ are known in some 
extent \cite{gra2, uga}. We also mention that there is an 
extensive literature about various submanifolds of $S^6$ equipped with 
various structures. Far from being complete, for instance the almost 
complex submanifolds of the Cayley almost complex $S^6$ are studied in 
\cite{gra1} and an excellent survey about the Lagrangian submanifolds of the 
Cayley nearly K\"ahler $S^6$ is \cite{chen}. 

The long resistance of the problem against proving non-existence may 
indicate that one should rather try to seek a complex structure on $S^6$.
However the irregular features of this hypothetical space $X$ convince 
us that {\it asymmetry, transcendental} 
(i.e., non-algebraic) {\it methods} and {\it inhomogeneity} should play a 
key role in finding it. 

Indeed, by a result of Wood \cite{woo} from a ``physical'' viewpoint the 
\ae sthetic Cayley almost complex structure is energetically remarkable 
unstable.\footnote{However for a contradictory result cf. \cite{bor-lam-sal}.} 
Therefore it is not surprising that if integrable almost complex 
structures on $S^6$ exist then they would appear ``far'' from the Cayley one. 
For instance Peng and Tang \cite{pen-tan} recently have constructed a novel 
almost complex structure by twisting and extending the standard complex 
structure on $\C^3\subset S^6$. It is a vast deformation of the Cayley 
almost complex structure and is integrable except a narrow equatorial 
belt in $S^6$. However it is still orthogonal with respect to the standard 
round metric hence LeBrun's theorem \cite{leb, zho1} forbids its full 
integrability.

One may try to seek $X$ by minimizing another energy functional over a 
space somehow related to $S^6$. Two questions arise: what this space should 
be and how the energy functional should look like?

As we have indicated above, the central difficulty with this hypothetical 
complex structure is that it should be compact and non-K\"ahler at the same 
time. However in fact plenty of such complex manifolds are known to exist for 
a long time: classical examples are the Hopf-manifolds \cite{hop} from 1948, 
the Calabi--Eckmann manifolds \cite{cal-eck} from 1953; or the compact, even 
dimensional, simply connected Lie groups---as it was observed by Samelson 
\cite{sam} also in 1953. {\it A fortiori} our candidate is the compact 
simply connected 14 dimensional exceptional Lie group ${\rm G}_2$ as a 
complex manifold which also arises as the total space of a non-trivial 
${\rm SU}(3)$ principal bundle over $S^6$. 

Regarding the second question a physisist is experienced that 
there is a close analogy between the existence of a {\it geometric structure} 
on a manifold on the mathematical side (cf. e.g. \cite{joy}) and the existence 
of a {\it spontaneous symmetry breaking} in classical Yang--Mills--Higgs-like 
systems on the physical side. This later process is familiar in both classical 
and quantum field theory \cite{che-li, kak, wei}. Therefore we are tempted to 
construct an appropriate (non-linear) field theory on ${\rm G}_2$ what 
we call a {\it Yang--Mills--Higgs--Nijenhuis theory}. Its dynamical 
variables are a Riemannian metric, a compatible gauge and a Higgs 
field with a usual Higgs potential but coupled to the gauge field through the 
Nijenhuis tensor; hence the Higgs field in the spontaneously broken classical 
vacuum state represents a complex structure on ${\rm G}_2$. Following the 
standard procedure used in dimensional reduction (cf. \cite{che-li, 
gre-sch-wit, kak, wei} or for a recent application \cite{wit1, wit2}) we 
Fourier expand this classical vacuum state with respect to ${\rm SU}(3)$ and 
expect that the ground Fourier mode---independent of the vertical 
directions---descends and gives rise to a classical vacuum solution 
in an effective Yang--Mills--Higgs--Nijehnuis field theory on $S^6$. If this 
classical vacuum solution exists then it would represent a complex structure 
on $S^6$. In checking that this is indeed the case we exploit the double role 
played by the group ${\rm SU}(3)$: it can be used not only for Fourier 
expansion but at the same time to put an almost complex structure onto 
$TS^6$ by constructing the tangent bundle as a vector bundle associated to 
the ${\rm SU}(3)$ principal bundle $\pi :{\rm G}_2\rightarrow S^6$ by the aid 
of the standard $3$ dimensional complex representation of ${\rm SU}(3)$. It 
also should become clear in the course of the construction why this approach 
breaks down in apparently similar situations like 
${\rm SO}(4n+1)/{\rm SO}(4n)\cong S^{4n}$. It turns out that, in sharp 
contrast to the case of ${\rm G}_2/{\rm SU}(3)\cong S^6$, in the former cases 
the tangent bundles $TS^{4n}$ cannot be constructed as {\it complex} 
associated vector bundles simply because the standard vector representations 
of ${\rm SO}(4n)$'s are {\it real}. Guided by these heuristic ideas in what 
follows we hope to prove rigorously the existence of a complex structure on 
the six-sphere. 

The paper is organized as follows. Sect. \ref{two} summarizes the physical 
background motivation. It contains a formulation of a hierarchy of classical 
field theories with gauge symmetry describing various spontaneous symmetry 
breakings related to the existence of various geometric structures on the 
underlying manifold. Apparently these sort of spontaneous symmetry breakings 
are not distinguished by physicists. In particular we identify a so-called 
``Yang--Mills--Higgs--Nijenhuis field theory'' on the tangent bundle of an 
even dimensional oriented Riemannian manifold describing a so-called ``weak 
spontaneous symmetry breaking'' which is the physical reformulation of the 
existence of a non-K\"ahlerian integrable almost complex structure on the 
underlying manifold. Unfortunately this theory is highly non-linear hence 
solving its field equations directly or even the detailed analytical study of 
them is far beyond our technical limits. 

Instead in Sect. \ref{three} we rigorously develop a sort of global Fourier 
expansion which makes it possible to push down sections of vector bundles 
with lifted $G$-action over the total space $P$ of a principal $G$-bundle to 
sections of quotient vector bundles over the base space $M=P/G$. However if we 
want to make a contact with the gauge symmetry inherently present in our 
formulation it turns out that this global Fourier expansion is not unique: 
Fourier expansion of a section is in general {\it not} gauge equivalent to the 
Fourier expansion of the gauge transformed section. Or saying equivalently, 
the Fourier expansion of vector-valued sections is not unique in general and 
depends on the choosen lifted $G$-action on the vector bundle on $P$. 

As a warming-up in Sect. \ref{five} we construct the standard 
complex structure on $S^2$ with our tools developed here. 
Actually no Fourier expansion is required for this trivial example. 

In Sect. \ref{four} then our tools are used for $S^6$ as follows: after 
identifying ${\rm G}_2$ with an ${\rm SU}(3)$ bundle over $S^6$ we Fourier 
expand the horizontal part of the underlying almost complex tensor field of a 
complex structure on ${\rm G}_2$. The ground Fourier mode depends on the 
chosen Fourier expansion and---since simply being the average along 
the fibers---always descends to $S^6$ as a skew-symmetric tensor field but in 
general not as an almost complex tensor field. However we demonstrate that 
there exists an ${\rm SU}(3)$ Fourier expansion which is truncated more 
precisely the corresponding ground Fourier mode coincides with the full 
horizontal component of the original integrable almost complex structure on 
${\rm G}_2$; hence it gives rise to an almost complex structure on $S^6$ 
(cf. Lemma \ref{vakuumlemma} here). We also make it clear why the
resulting almost complex structure on $S^6$ cannot be homogeneous in spite of 
the fact that the complex structure it stems from is homogeneous on 
${\rm G}_2$. In short the reason is that the tangent spaces of $S^6$ are 
identified with the horizontal part of $T{\rm G}_2$ not simply by the 
derivative of the projecion $\pi :{\rm G}_2\rightarrow S^6$ as usual but 
instead in a fiberwise twisted manner (cf. Lemma \ref{homogenlemma} here). 
This section contains our main result, namely that this non-homogeneous almost 
complex structure is integrable hence there exists a complex manifold 
homeomorphic to the six-sphere (cf. Theorem \ref{Xtetel} here). 


\section{Geometric structures and symmetry breaking}
\label{two}


In this section we construct classical field theories with 
gauge symmetry whose spontaneously broken vacua give rise to important 
geometric structures on the underlying manifold. In particular a so-called 
Yang--Mills--Higgs--Nijenhuis field theory is exhibited with a spontaneously 
broken vacuum corresponding to a non-K\"ahlerian complex structure. 

Let $M$ be an $n$ dimensional real manifold and pick an abstract Lie group 
$G\subset{\rm GL}(n,\R )$. Recall that \cite[Section 2.6]{joy} a {\it 
$G$-structure} on $M$ is a principal sub-bundle of the frame bundle of 
$M$ whose fibers are $G$ and the structure group of this sub-bundle is also 
$G$. A $G$-structure is called {\it integrable} (or {\it torsion-free}) if 
there is a torsion-free connection $\nabla$ on $TM$ with ${\rm Hol}\nabla$ 
being (a subgroup of) $G$. An integrable $G$-structure is interesting because 
the individual infinitesimal geometric structures on the tangent spaces 
stem from an underlying global non-linear structure on $M$.

As a starting setup suppose that a $G$-structure on $M$ is given 
together with some connection on $TM$ whose holonomy group satisfies 
${\rm Hol}\nabla \subseteqq G$. Note that $G$ plays a 
double role here. Pick a Lie subgroup $H\subset G$. We are going 
to describe two broad versions of a ``spontaneous symmetry breaking'' over 
$M$ as follows. Our first concept is a {\it strong spontaneous symmetry 
breaking}, denoted by $G\Rrightarrow H$, in which {\it $G$ 
as the holonomy group of $\nabla$} is reducible to (a subgroup of) $H$. 
That is, the connection satisfies that ${\rm Hol}\nabla\subseteqq H$.

To clarify this concept, an example can be constructed as follows. Take 
a $2m$ real dimensional oriented Riemannian manifold $(M,g)$. 
Then $TM$ is an ${\rm SO}(2m)$ vector bundle associated with the 
standard representation of ${\rm SO}(2m)$ and has a gauge group $\cg_{TM}$ 
consisting of fiber-preserving isomorphisms of $TM$ which keep 
orientation and are orthogonal with respect to the metric $g$ along each fiber.
Let $\nabla$ be an arbitrary SO$(2m)$ 
connection on $TM$ and let $R_\nabla$ and $T_\nabla$ be its curvature and 
torsion respectively.

Consider moreover another associated bundle to the frame bundle 
making use of the adjoint representation of SO$(2m)$ on its Lie algebra. 
The fibers of this bundle ${\rm Ad}M\subset{\rm End}TM$ are 
${\mathfrak s}{\mathfrak o}(2m)\subset {\rm End}\:\R^{2m}$. Let 
$\Phi\in C^\infty (M;{\rm Ad}M)$ be its generic section. A geometric 
example of a section of this kind is an orthogonal almost complex 
structure whose induced orientation agrees with the given one. Indeed, 
since for a real oriented vector space with a scalar product 
$(V, \langle\cdot\:,\:\cdot\rangle )$ the moduli space of orthogonal 
complex structures $J$ compatible with the orienation is precisely ${\rm 
SO}(V)\cap{\mathfrak s}{\mathfrak o}(V)\subset{\rm End}\:V$, it follows that 
we can look at any orthogonal almost complex structure with respect to $g$ 
and compatible induced orientation as a section 
$J=\Phi\in C^\infty (M;{\rm Ad}M)$. Let 
$\nabla :C^\infty (M;{\rm Ad}M)\rightarrow 
C^\infty (M;{\rm Ad M}\otimes\wedge^1M)$ be the connection on ${\rm Ad}M$, 
associated to the SO$(2m)$ connection on $TM$ introduced above. Also fix a 
real number $e\not= 0$.

Then for instance an ${\rm SO}(2m)\Rrightarrow{\rm U}(m)$ strong spontaneous 
symmetry breaking arises by seeking the hypothetical vacuum solutions 
(see below) of the energy functional
\[\ce (\nabla ,\Phi ):=
\frac{1}{2}\int\limits_M\left(\frac{1}{e^2}\vert 
R_\nabla\vert_g^2+\vert\nabla\Phi\vert_g^2+e^2\vert\Phi\Phi^*-{\rm 
Id}_{TM}\vert^2_g\right)\dd V.\]
Various pointwise norms induced by the metric on various  SO$(2m)$ 
vector bundles and the Killing form on
${\mathfrak s}{\mathfrak o}(2m)$ are used here to calculate norms of 
tensor fields in the integral. We note that $\Phi^*$ is the pointwise 
adjoint (i.e., transpose) of the bundle map $\Phi :TM\rightarrow TM$ with 
respect to $g$. 

$\ce (\nabla ,\Phi )\geqq 0$ and it is 
invariant under ${\rm SO}(2m)$ acting as gauge transformations 
$\alpha\in\cg_{TM}$ 
given by $(\nabla ,\Phi )\mapsto (\alpha\nabla\alpha^{-1} ,
\alpha\Phi\alpha^{-1})$. In fact $\ce (\nabla ,\Phi )$ is the usual 
Euclidean action of a Yang--Mills--Higgs theory from the physics 
literature formulated on $TM$. Hence, using 
physicists' terminology \cite{che-li,kak,wei} the above functional can be 
regarded as the action functional of a {\it spontaneously broken $2m$
dimensional Euclidean ${\rm SO}(2m)$ Yang--Mills--Higgs system} on $TM$.
In this context $\nabla$ is the ``${\rm SO}(2m)$ gauge 
field'' and $\Phi$ is the ``Higgs field in the adjoint representation'', 
$e$ represents the ``coupling constant'' and the quartic polynomial 
$\vert\Phi\Phi^*-{\rm Id}_{TM}\vert^2_g$ is the ``Higgs potential''. 
If a smooth vacuum solution $\ce (\nabla ,\Phi ,g)=0$ exists then 
\[\left\{\begin{array}{ll}
0= & R_{\nabla}\\
0= & \nabla\Phi\\
0= & \Phi\Phi^*-{\rm Id}_{TM}.
        \end{array}\right.\]
(Note that these are {\it not} the Euler--Lagrange equations of the 
system!) The first equation says that $\nabla$ is flat. The third 
equation dictates $\Phi^2=-{\rm Id}_{TM}$ hence restricts the Higgs 
field to be an orthogonal almost complex structure $\Phi =J$ with 
respect to $g$. The second equation reduces the holonomy of 
$\nabla$ in the usual way \cite[Proposition 2.5.2]{joy} to 
(by flatness necessarily a discrete subgroup of) ${\rm U}(m)$. 
Therefore this vacuum---if exists---describes a (not necessarily integrable) 
almost complex structure on $M$. The fact that for this vacuum solution 
${\rm Hol}\nabla\subset{\rm U}(m)\subset{\rm SO}(2m)$ justifies 
that this is indeed an example for a strong spontaneous symmetry breaking 
${\rm SO}(2m)\Rrightarrow{\rm U}(m)$. The resulting structure is 
non-integrable and flat hence mathematically less interesting. We have 
included this example because of its relevance in physics.

This ``illness'' of the vacuum is cured if the curvature is 
replaced by the torsion. In passing consider a new theory with
\[\ce (\nabla ,\Phi ,g):=\frac{1}{2}\int\limits_M\left(\frac{1}{e^2}
\vert {\mathfrak S}(R_\nabla)\vert^2_g+
\vert\nabla g\vert^2_g+\vert\nabla\Phi\vert_g^2+
e^2\vert\Phi\Phi^*-{\rm Id}_{TM}\vert^2_g\right)\dd V\]
where ${\mathfrak S}(R_\nabla )$ is the symmetrization of $R_\nabla$ given by
\[{\mathfrak S}(R_\nabla )(u,v)w:={\mathfrak S}(R_\nabla (u,v)w)=R_\nabla 
(u,v)w+R_\nabla (w,u)v+R_\nabla (v,w)u\]
for vector fields $u,v,w$ on $M$. Bianchi's first identity says that
\[{\mathfrak S}(R_\nabla (u,v)w)={\mathfrak S}(T_\nabla (T_\nabla 
(u,v),w)+\nabla_uT_\nabla (v,w))\]
demonstrating that ${\mathfrak S}(R_\nabla )=0$ if and only if $T_\nabla =0$. 
This functional is again non-negative and taking into account that 
$\alpha^*g=g$ it is acted upon by $\cg_{TM}$ as $(\nabla ,\Phi ,g)\mapsto 
(\alpha\nabla\alpha^{-1} ,\alpha\Phi\alpha^{-1} ,g)$ 
leaving the functional invariant. Although this theory looks physically less 
obvious this departure from physics yields a more familiar structure on the 
mathematical side. Namely, a hypothetical smooth minimizing 
configuration $\ce (\nabla ,\Phi ,g)=0$ gives rise to a K\"ahler structure on 
$M$. Note that ${\rm Hol}\nabla\subseteqq{\rm U}(m)\subset{\rm SO}(2m)$ hence 
this is again an (integrable) strong spontaneous symmetry breaking 
${\rm SO}(2m)\Rrightarrow{\rm U}(m)$.

Returning to the original setup, our second concept is a {\it weak 
spontaneous symmetry breaking}, denoted as $G\Rightarrow H$, in 
which $G$ {\it as the structure group of $TM$} is reducible to (a 
subgroup of) $H$. That is, simply there is a $H$-structure on $M$ but probably 
the connection is not compatible with this structure. A strong 
symmetry breaking always implies the corresponding 
weak one cf. e.g. \cite[Theorem 2.3.6]{joy} however not the other way round. 

An (integrable) weak symmetry breaking ${\rm SO}(2m)\Rightarrow{\rm U}(m)$ 
is provided by the following example which goes along the lines of the 
previous cases. Consider the first order quadratic differential operator
\[N_\nabla : C^\infty (M;{\rm Ad}M)\longrightarrow
C^\infty (M; ({\rm End}TM)\otimes\wedge^1M)\cong C^\infty (M; 
TM\otimes\wedge^1M\otimes\wedge^1M)\]
whose shape on two vector fields $u,v$ on $M$ is defined to be
\[(N_\nabla\Phi )(u,v):=(\nabla_{\Phi u}\Phi )v-\Phi (\nabla_u\Phi )v-
(\nabla_{\Phi v}\Phi )u +\Phi (\nabla_v\Phi )u\]
or equivalently in the particular gauge given by a local coordinate system 
$(U, x^1,\dots,x^{2m})$ it looks like \cite[Eqn. 15.2.5]{gre-sch-wit}
\begin{equation}
\left(\left(N_\nabla\Phi\right)\vert_U\right)^k_{\:\:\:ij}=
\sum\limits_{l=1}^{2m}
\Phi^l_{\:\:i}(\nabla_l\Phi^k_{\:\:j}-\nabla_j\Phi^k_{\:\:l})-
\sum\limits_{l=1}^{2m}\Phi^l_{\:\:j}(\nabla_l\Phi^k_{\:\:i}-  
\nabla_i\Phi^k_{\:\:l}),\:\:\:\:\:i,j,k=1,\dots,2m.
\label{lokalis-nijenhuis}
\end{equation}
Note that $N_\nabla\Phi$ transforms as a $(2,1)$-tensor under ${\rm SO}(2m)$ 
acting on $TM$. This time set
\begin{equation} 
\ce (\nabla ,\Phi ,g):=\frac{1}{2}\int\limits_M\left(\frac{1}{e^2}
\vert {\mathfrak S}(R_\nabla )\vert^2_g+
\vert\nabla g\vert^2_g+\vert N_\nabla\Phi\vert_g^2+
e^2\vert\Phi\Phi^*-{\rm Id}_{TM}\vert^2_g\right)\dd V.
\label{YMH}
\end{equation}
This functional is again invariant under $\cg_{TM}$ acting on 
$(\nabla ,\Phi ,g)$ as before and it is non-negative. Therefore it describes 
a sort of gauge theory on $M$ what we call a 
{\it Yang--Mills--Higgs--Nijenhuis field theory}. The peculiarities here 
compared to the usual Yang--Mills--Higgs Lagrangian are that the curvature 
has again been replaced with the torsion, the 
connection is minimally coupled to the metric but its usual minimal  
coupling $\nabla\Phi$ to the Higgs field has been replaced with $N_\nabla\Phi$ 
regarded as a non-minimal quadratic coupling in (\ref{YMH}). This further 
departure from the physical side gives rise again to a familiar structure on 
$M$ on the mathematical side as follows.

If the smooth vacuum $\ce (\Phi ,\nabla ,g)=0$ is achieved in (\ref{YMH}) 
over $M$ then it satisfies
\begin{equation}
\left\{\begin{array}{ll}
0= & {\mathfrak S}(R_\nabla )\\
0= & \nabla g\\
0= & N_\nabla\Phi\\
0= & \Phi\Phi^*-{\rm Id}_{TM}.
        \end{array}\right.
\label{vakuum}
\end{equation}
The first and second equations say that the ${\rm SO}(2m)$ connection 
$\nabla$ is the Levi--Civita connection of $g$. The fourth equation 
ensures us that $\Phi^2=-{\rm Id}_{TM}$ i.e., the Higgs field reduces to an 
orthogonal almost complex structure $\Phi =J$ with respect to the metric and 
is compatible with the orientation. Exploiting this together with the fact 
that the Levi--Civita connection is torsion-free the operator 
$N_\nabla\Phi =N_\nabla J$ in the gauge (\ref{lokalis-nijenhuis}) cuts down to
\[\left(\left(N_{\nabla}J\right)\vert_U\right)^k_{\:\:\:ij}=
\sum\limits_{l=1}^{2m}J^l_{\:\:i}\left(\partial_lJ^k_{\:\:j}-
\partial_jJ^k_{\:\:l}\right) 
-\sum\limits_{l=1}^{2m}J^l_{\:\:j}\left(\partial_lJ^k_{\:\:i}-
\partial_iJ^k_{\:\:l}\right),\:\:\:\:\:i,j,k=1,\dots,2m\]
hence $N_{\nabla} J=N_J$ is the {\it Nijenhuis tensor} of the almost complex 
structure $J$ given by
\[N_J(u,v)= -[u,v]+[Ju, Jv]-J[u, Jv]-J[Ju, v].\]
Consequently the third equation of (\ref{vakuum}) says that $N_J=0$ that is, 
this almost complex structure is integrable on $M$ and is orthogonal with 
respect to $g$. Therefore if a vacuum solution to (\ref{YMH}) exists it makes 
$M$ into a complex analytic $m$-space in light of the Newlander--Nirenberg 
theorem \cite{new-nir}. Note that this time the structure group of $TM$ 
has been cut down to (a subgroup of) ${\rm U}(m)$ but probably for the 
connection ${\rm Hol}\nabla\not\subseteqq{\rm U}(m)$ hence we indeed 
obtain an example for an (integrable) weak spontaneous symmetry breaking 
${\rm SO}(2m)\Rightarrow {\rm U}(m)$. Accordingly the resulting complex 
structure on $M$ is more general than the K\"ahler one in the previous example.

We have described two different (i.e., strong and weak) symmetry 
breakings from ${\rm SO}(2m)$ down to ${\rm U}(m)$ over a Riemannian manifold. 
However for physicists these concepts coincide because for them a ``spontaneous 
symmetry breaking'' has a slightly different meaning. It simply means that 
although one starts with an ${\rm SO}(2m)$-invariant theory over $M$, 
the vacuum $(\nabla ,\Phi ,g)$ is stabilized by a subgroup of ${\rm SO}(2m)$ 
isomorphic to ${\rm U}(m)$ only. Of course this criteria holds 
in all of our cases. 

We know that all integrable solutions to the previous strong 
problem---i.e., the K\"ahler structures---provide integrable solutions 
to the weak one---i.e., complex analytic structures---too. But we 
rather raise the question whether or not there are solutions strictly to 
the weaker problem which do not stem from strong solutions. The answer is 
yes, even in the compact case. In what follows let $G$ be a connected compact 
even dimensional Lie group. Due to Samelson \cite{sam} homogeneous complex 
structures on $G$ can be constructed as follows. Write ${\mathfrak g}:=T_eG$ 
for the Lie algebra of $G$ and let ${\mathfrak g}^\C$ be its complexification. 
A complex Lie subalgebra ${\mathfrak s}\subset{\mathfrak g}^\C$ is called a 
{\it Samelson subalgebra} if 
$\dim_\C{\mathfrak s}= \frac{1}{2}\dim_\C({\mathfrak g}^\C )$ and 
${\mathfrak s}\cap{\mathfrak g}=0$ where now ${\mathfrak g}
\subset{\mathfrak g}^\C$. Infinitesimally ${\mathfrak s}$ 
determines an almost complex structure 
$J_e:{\mathfrak g}\rightarrow {\mathfrak g}$ which satisfies that 
${\mathfrak s}$ is the $-\ii$-eigenspace of its complex linear extension 
$J_e^\C :{\mathfrak g}^\C\rightarrow {\mathfrak g}^\C$. Globally 
${\mathfrak s}\subset{\mathfrak g}^\C$ gives rise to a complex 
Lie subgroup $S\subset G^\C$ of the complexification and there exists 
a diffeomorphism $G^\C\cong G\times S$; a complex structure on $G$ 
induced by ${\mathfrak s}$ arises by taking $G\cong G^\C /S$. 
Left translations act by biholomorphisms hence $G$ is a compact 
homogeneous complex manifold. If $H^2(G;\C )\cong 0$ the space $G$ 
{\it cannot} be projective algebraic (cf. e.g. \cite[Lemma 4.1]{cam-dem-pet}) 
or even K\"ahler hence $\nabla J\not =0$ with respect to any metric. An 
example is ${\rm SU}(2)\times{\rm SU}(2)$ as a complex 
$3$-manifold.\footnote{Note that the smooth manifold $S^3\times S^3$ 
can be given the structure of a complex manifold in the sense of 
Calabi--Eckmann \cite{cal-eck}, too. Although in principle this complex 
structure differs from the one stemming from the Samelson construction 
presented here they are in fact isomorphic \cite{dau,mag}.} Another 
example is the exceptional Lie group ${\rm G}_2$ as a complex $7$-manifold.

As a $T$-module with respect to the adjoint action of a maximal torus 
$T\subseteqq G$ we have the usual decomposition
\begin{equation}
{\mathfrak g}={\mathfrak t}\oplus\bigoplus\limits_{\mu\in R^+}
{\mathfrak g}_\mu
\label{gyokok}
\end{equation}
where ${\mathfrak t}={\mathfrak t}^\C\cap{\mathfrak g}$ and 
${\mathfrak g}_\mu =({\mathfrak g}^\C_\mu\oplus 
{\mathfrak g}^\C_{-\mu})\cap{\mathfrak g} =
{\mathfrak g}_{-\mu}$ and $R^+$ is 
supposed to contain exactly one element from each pair $\{\mu ,-\mu\}$ 
of real roots. The almost complex operator 
$J_e$ has a corresponding decomposition. As a real vector space 
${\mathfrak t}\cong\R^l$ where $l=\dim_\R{\mathfrak t}$ 
is the rank of $G$. Since $\dim_\R G$ is even, $l$ is even hence $J_e$ 
restricted to ${\mathfrak t}$ gives a complex structure $J_0$ on 
${\mathfrak t}$ providing an isomorphism ${\mathfrak t}\cong\C^{\frac{l}{2}}$. 
In addition, on any ${\mathfrak g}_\mu\cong\R^2$ we obtain a 
unique complex structure $J_\mu$ on $\R^2$ providing 
${\mathfrak g}_\mu\cong\C$. Consequently at the unit $e\in G$ (hence 
everywhere) the integrable almost complex operator $J$ representing the 
complex structure splits:
\begin{equation}
J_e=J_0\oplus\bigoplus\limits_{\mu\in R^+}J_\mu .
\label{J-szetszedes}
\end{equation}

Suppose $G$ is an even dimensional compact Lie group with $H^2(G;\C )\cong 0$ 
and take $M:=G$ to be the base manifold in (\ref{YMH}). Picking a 
representative $J$ of the complex structure on $M$ put an induced orientation 
as well as a metric $g_J$ onto $M$ for which $J$ is orthogonal. Such metric 
exists by averaging any metric $g$ with respect to $J$ i.e., setting 
$g_J(u,v):=\frac{1}{2}g(u,v)+\frac{1}{2}g(Ju,Jv)$ for vector fields $u,v$ on 
$M$. Consequently taking $\Phi :=J$ and $g:=g_J$ and $\nabla$ to be its 
Levi--Civita connection we obtain that $(\nabla ,\Phi ,g)$ solves 
(\ref{vakuum}) on the compact space $M$ i.e., it gives rise to a strictly 
weak spontaneous symmetry breaking ${\rm SO}(2m)\Rightarrow{\rm U}(m)$ where 
$\dim_\R M=2m$. The spontaneously broken vacuum represents the complex 
structure on $M$ which is not K\"ahler.

Finally we clarify a subtle conceptual ambiguity: how should one interpret 
the effect of a general gauge transformation on a given geometric structure? 
Consider $(\R^{2m}, J)\cong\C^m$ and pick an element $A\in{\rm GL}(2m,\R )$. 
On the one hand if $A$ acts directly on $\R^{2m}$ by 
$(x^1,\dots, x^{2m})\mapsto A(x^1,\dots,x^{2m})$ then in this new frame $J$ 
looks like $AJA^{-1}$. In this case $A$ just represents a linear 
coordinate transformation on $\R^{2m}$ hence ``$AJA^{-1}$ is the same old 
complex structure in a new frame''. Globally, $f\in{\rm Diff}(M)$ 
acts on a vector field by the chain rule as $u\mapsto f_*u$ hence on $J$ by 
$J\mapsto f_*Jf^{-1}_*$ regarded as a ``coordinate transformation''. 
Consequently the transformed structure is considered being equivalent to $J$. 
More generally an $\alpha\in\cg_{TM}$ from the gauge group acts the same way: 
$J\mapsto \alpha J\alpha^{-1}$ therefore we obtain an embedding 
${\rm Diff}(M)\subset
{\rm Aut}(C^\infty (M;{\rm End}TM))\cong\cg_{TM}/\{\pm {\rm Id}_{TM}\}$. 
This natural embedding and the 
gauge theoretic formulation developed here dictates to extend this picture 
from diffeomorphisms to general gauge transformations as well. On the other 
hand if $A$ does not act on $\R^{2m}$ but acts on $J$ itself by 
$J\mapsto AJA^{-1}$ then $A$ describes a deformation of the complex structure 
hence ``$AJA^{-1}$ is a new complex structure in the same old frame''. 
Correspondingly, in this picture $\cg_{TM}$ acts by global deformations 
$J\mapsto \alpha J\alpha^{-1}$ of the almost complex structure. 
\vspace{0.1in}

\noindent {\bf Principle.} {\it Over an almost complex manifold there 
exists an abstract group operating in two essentially different ways. 
First, this group as the gauge group $\cg_{TM}$ describes 
the symmetries of a fixed almost complex manifold $(M,J)$ by gauge 
transforming everything. Second, this group as a "deformation group" can 
also describe deformations of a given almost complex manifold 
$(M,J)$ by deforming it into a new one $(M,J')$ and leaving other objects
intact.}
\vspace{0.1in}

\noindent Therefore if this group acts as the gauge group describing the 
symmetries of a fixed almost complex manifold $(M,J)$ by gauge transforming 
everything then in particular the Nijenhuis tensor 
(\ref{lokalis-nijenhuis}) transforms as a $(2,1)$-tensor under such a 
symmetry transformation:
\[N_J=N_\nabla J\longmapsto N_{\alpha\nabla\alpha^{-1}}(\alpha 
J\alpha^{-1} )=\alpha (N_\nabla J)=\alpha (N_J)\]
where $\alpha (T)$ denotes the natural action of $\alpha\in\cg_{TM}$ on a 
$(k,l)$-type tensor $T$. But if this group acts as a "deformation 
group" by deforming $J$ alone and not changing anything else then 
in particular the Nijenhuis tensor does not behave tensorially under such 
deformations:
\[N_J = N_\nabla J\longmapsto N_\nabla (\alpha J\alpha^{-1})=
N_{\alpha J\alpha^{-1}}\not=\alpha (N_J).\]
These two roles played by the gauge group should not be mixed in the 
forthcoming paragraphs.


\section{Fourier series expansion over principal bundles}
\label{three}


In this technical section we we collect some 
useful facts about expansion of functions into Fourier series over a 
connected compact Lie group $G$. A standard reference is for instance 
\cite{hel, may}. Then we generalize this to Fourier expand sections of vector 
bundles admitting lifted $G$-actions over principal $G$-bundles. These results 
also exist in the literature although in a somewhat implicit form (cf. 
e.g. \cite[Chapters II and III]{bro-tom}). This generalized Fourier 
expansion is a rigorous mathematical tool for performing dimensional 
reduction in classical and quantum field theories (cf. e.g. \cite{wit1, wit2}).

Let $G$ be a compact Lie group and let us denote by ${\rm Irr}(G;\C )$ the 
set of isomorphism classes of its finite dimensional complex irreducible 
representations. That is, if $\rho\in {\rm Irr}(G;\C )$ then there exists 
a finite dimensional complex vector space $V_\rho$ and a homomorphism 
$\rho : G\rightarrow{\rm Aut}V_\rho$. For a fixed (isomorphism 
class of) $\rho$ pick a basis in $V_\rho$ to write $(V_\rho 
,e_1,\dots,e_{\dim_\C V_\rho})\cong\C^{\dim_\C V_\rho}$ and denote the 
corresponding matrix elements of $\gamma\in G$ as $\rho_{ij}(\gamma )$. 
Then $\rho_{ij}:G\rightarrow\C$ is a continuous function and the 
Peter--Weyl theorem \cite[Theorem III.3.1]{bro-tom} implies that
\[\left\{\sqrt{\dim_\C V_\rho}\:\rho_{ij}:G\rightarrow\C\:\vert\:\rho\in 
{\rm Irr}(G;\C );\:\: i,j=1,\dots,\dim_\C V_\rho\right\}\]
forms a complete set of orthonormal functions in the Hilbert space 
$L^2(G;\C )$ of square-integrable complex-valued functions on $G$ with 
the usual scalar product 
$(f_1,f_2)_{L^2(G)}:=({\rm Vol}\:G)^{-1}\int_Gf_1\overline{f}_2\dd\gamma$. 
Here $\dd\gamma$ is a bi-invriant Haar measure on $G$ and 
${\rm Vol}\:G=\int_G\dd\gamma$. A function $f\in L^2(G;\C )$ then can be 
written as a formal sum
\[f(\gamma )\sim\sum\limits_{\rho\in {\rm Irr}(G;\C )}
\sqrt{\dim_\C V_\rho}\sum\limits_{i,j=1}^{\dim_\C V_\rho}a_{\rho ,ij}\:
\rho_{ij}(\gamma )\]
where the constants $a_{\rho ,ij}\in\C$ are given by  
\[a_{\rho ,ij}:=\frac{1}{{\rm Vol}\:G}\:\int\limits_G 
f(\gamma )\:\overline{\rho}_{ij}(\gamma )\dd\gamma .\]

This expansion can be obtained in a more invariant (and abstract) way by 
the aid of the irreducible representations themselves without referring to the 
individual matrix elements. Right translation on $G$ induces an 
infinite dimensional unitary representation of $G$ from the left on 
$L^2(G;\C )$ given by 
\[(\gamma\cdot f)(\delta ) :=f(\delta\gamma )\] 
for all $\delta\in G$. This representation---the {\it regular 
representation}---gives rise to an orthogonal decomposition
\[L^2(G;\C )\cong\bigoplus\limits_{\rho\in{\rm Irr}(G;\C )}W (V_\rho )\]
where $W (V_\rho )$ is the isotypical summand for the finite dimensional 
irreducible representation $V_\rho$ i.e., $W (V_\rho )$ is 
the $L^2$-closure of the span of all subspaces in $L^2(G;\C )$ isomorphic 
to $V_\rho$ as $G$-modules (in fact in the case at hand $W(V_\rho )$ is 
the direct sum of $\dim_\C V_\rho$ copies of $V_\rho$). 
The continuous orthogonal projection $\Pi_\rho :L^2(G;\C )\rightarrow W 
(V_\rho )$ is given by the character
\[\chi_\rho :={\rm Tr}\:\rho =\sum_{i=1}^{\dim_\C V_\rho}\rho_{ii}\]
of the corresponding representation as follows:
\begin{equation}
\Pi_\rho f :=\frac{\dim_\C V_\rho}{{\rm Vol}\:G}\int\limits_G(\gamma\cdot f)\:
\overline{\chi}_\rho (\gamma )\dd\gamma .
\label{egyutthatok}
\end{equation}  
The associated formal sum 
\begin{equation}
f \sim\sum_{\rho\in{\rm Irr}(G;\C )}\Pi_\rho f
\label{fourier}
\end{equation}
is called the {\it Fourier expansion of $f$ with respect to $G$}. If $f$ 
is smooth then $\Pi_\rho f$'s are also smooth and (\ref{fourier}) converges 
uniformly and pointwise equality holds.

Now we proceed further and construct more general Fourier expansions. Let $M$ 
be an arbitrary connected smooth manifold. As usual let $G$ be a 
connected compact Lie group with a fixed bi-invariant Haar measure $\dd\gamma$ 
and consider a principal $G$ bundle $\pi : P\rightarrow M$. This means that 
$G$ has a free right action on $P$ such that $P/G\cong M$. Suppose that there 
is a Riemannian metric on $P$. Also let $E$ be a complex vector bundle 
$p:E\rightarrow P$ such 
that the right action of $G$ on $P$ lifts to $E$ rendering it a vector bundle 
with some fixed smooth right $G$-action. Put some $G$-invariant Hermitian 
scalar product onto $E$; this together with the metric on $P$ 
gives a complex Hilbert space $L^2(P;E)$. A generic section is denoted by 
$s\in L^2(P;E)$. The right action of $G$ on $P$ and its lifted right 
action on $E$ induces a continuous representation of $G$ on 
$C^\infty (P;E)$ from the left given by 
\[(\gamma\cdot s)(y):=s(y\gamma )\gamma^{-1}\] 
whose unique continuous extension makes $L^2(P;E)$ into an infinite 
dimensional complex unitary $G$-module. We get again an orthogonal 
decomposition into isotypical summands
\[L^2(P;E)\cong\bigoplus\limits_{\rho\in{\rm Irr}(G;\C )}W(V_\rho )\]
as before with continuous orthogonal projections 
$\Pi_\rho :L^2(P; E)\rightarrow 
W(V_\rho )$ given by fiberwise integration\footnote{Although there are 
no canonical isomorphisms between the fibers of $\pi :P\rightarrow M$ and $G$ 
the Haar measure on $G$ can be pulled back to the fibers with any of these 
isomorphisms in an unambigous way taking into account its translation 
invariance. Hence we obtain well-defined measures on the fibers.}
\begin{equation}
\Pi_\rho s :=\frac{\dim_\C V_\rho}{{\rm Vol}\:G}\int\limits_G(\gamma\cdot s)\:
\overline{\chi}_\rho (\gamma )\dd\gamma .
\label{projekcio}
\end{equation}
These considerations suggest to define the {\it global Fourier expansion with 
respect to $G$} of a section $s\in L^2(P;E)$ as the formal sum
\begin{equation}
s\sim\sum\limits_{\rho\in {\rm Irr}(G;\C )}\Pi_\rho s.
\label{globalis.fi.fourier}
\end{equation}
If it happens that $M$ is compact and $s\in C^\infty (P;E)$ then it 
follows from the general theory of Fourier expansions \cite[Theorems 7.1 or 
8.6]{may} that also $\Pi_\rho s\in C^\infty (P; E)$ for all 
representations and (\ref{globalis.fi.fourier}) converges uniformly over $P$ 
and pointwise equality holds.

The ground Fourier mode $s_1:=\Pi_1s\in L^2(P;E)$ corresponding to the
trivial representation $\rho\cong 1$ of $G$ is nothing else than the 
average of $s$ along the fibers. More precisely exploiting left-translation 
invariance of the Haar measure it satisfies
\begin{eqnarray}
s_1(y\delta )&=&\frac{1}{{\rm Vol}\:G}\int\limits_G\left(\gamma\cdot 
s\right)(y\delta)\dd\gamma\nonumber\\
&=&\frac{1}{{\rm Vol}\:G}\int\limits_Gs(y\delta\gamma )\gamma^{-1}\dd\gamma
=\frac{1}{{\rm Vol}\:G}\int\limits_Gs(y\gamma )(\delta^{-1}\gamma )^{-1}
\dd\gamma =
\left(\frac{1}{{\rm Vol}\:G}\int\limits_Gs(y\gamma )\gamma^{-1}
\dd\gamma\right)\delta\nonumber\\
& =&s_1(y)\delta\nonumber
\end{eqnarray}
for all $\delta\in G$. Making use of the right $G$-action on $E$ we 
can form the natural collapsed bundle $E/G$ over $P/G=M$. A vector in the 
fiber $(E/G)_x$ over $x\in M$ corresponds to a section of 
$E\vert_{\pi^{-1}(x)}$ consisting of the equivalence class of vectors with 
respect to the right $G$-action. Hence $s_1$ descends uniquely to a section 
$\tilde{s}_1: M\rightarrow E/G$ of the collapsed bundle satisfying 
$\pi^*\tilde{s}_1 =s_1$. 

By the aid of the Hermitian structure on $E$ we have a gauge group $\cg_E$. 
If $\alpha\in\cg_E$ then it acts on sections by $s\mapsto \alpha s$ as usual. 
However we have seen that when performing Fourier expansions these sections are 
also acted upon by the group $G$ as $s\mapsto \gamma\cdot s$ constructed 
above. This motivates to let $\cg_E$ act on $G$-actions as 
$\gamma\mapsto\gamma_\alpha$ given by
\[(\gamma_\alpha\cdot s)(y):=s(y\gamma )\gamma_\alpha^{-1}\]
where the $\alpha$-twisted lifted right $G$-action 
$\gamma_\alpha :E_y\rightarrow E_{y\gamma}$ has the form
\[s(y)\gamma_\alpha :=\alpha (y\gamma )\left( (\alpha^{-1}(y)s(y))\gamma 
\right) .\]
This yields an identity
\[\gamma_\alpha\cdot (\alpha s)=\alpha (\gamma\cdot s).\]
With respect to the $\alpha$-twisted lifted right $G$-action on $E$ we can 
form again the collapsed bundle what we denote by $E/_\alpha G$. 
Since by assumption the Hermitian structure on $E$ is $G$-invariant we also 
obtain a gauge group $\cg_{E/_\alpha G}$ yielding a subgroup 
\[\cg_G:=\pi^*\cg_{E/_\alpha G}\subset\cg_E\] 
consisting of gauge transformations which are constant along the 
fibers.\footnote{Note that $\cg_G\subset\cg_E$ is independent of 
$\alpha\in\cg_E$.} 
Consequently if $\Pi_{\alpha ,\rho} (\alpha s)$ denotes a Fourier mode of 
$\alpha s$ with respect to $s\mapsto \gamma_\alpha\cdot s$ and 
$\Pi_{{\rm Id}_E,\rho} s:=\Pi_\rho s$ is that of $s$ with respect to $s\mapsto 
\gamma\cdot s$ then we can see from (\ref{projekcio}) that if in particular 
$\alpha\in\cg_G$ then $\Pi_{\alpha ,\rho} (\alpha s)=
\alpha (\Pi_{{\rm Id}_E,\rho} s)$. This means that the two 
Fourier expansions (\ref{globalis.fi.fourier}) are compatible with the gauge 
transformation. Since the twisted $G$-action $s\mapsto \gamma_\alpha\cdot s$ 
is a $G$-action on $L^2(P; E)$ on its own right we make the following 

\begin{definition}
Pick two gauge transformations $\alpha ', \alpha ''\in\cg_E$ and consider two 
$G$-actions on $L^2(P;E)$ of the form $s\mapsto \gamma_{\alpha '}\cdot s$ and 
$s\mapsto \gamma_{\alpha ''}\cdot s$ respectively. These give 
rise to two global Fourier expansions (\ref{globalis.fi.fourier}) with 
respect to $G$ on $L^2(P;E)$ given by the straightforwardly modified 
projections (\ref{projekcio}) respectively. 

These Fourier expansions are called {\em gauge equivalent} if there 
exists a $\beta\in \cg_G\subset\cg_E$ satisfying $\alpha '' =\beta\alpha '$.
\label{ekv.fourier}
\end{definition}

\begin{remark}\rm 1. For different $\alpha\in\cg_E$ the $G$-module
structures on $L^2(P;E)$ are unitarily equivalent. However the induced 
Fourier expansions are not equivalent and their moduli space 
is the quotient $\cg_G\backslash\cg_E$ in some sense. Also notice that two 
gauge equivalent Fourier expansions on $E$ over $P$ give rise to gauge 
equivalent ground modes on $E/_{\alpha '}G\cong E/_{\alpha ''}G$ over 
$P/G=M$ with respect to $\cg_{E/_{\alpha '}G}\cong\cg_{E/_{\alpha ''}G}$ 
which is the collapsed gauge group. Hence taking the ground mode is meaningful. 

To simplify notation given $\alpha\in\cg_E$ and its equivalence class 
$[\alpha ]\in \cg_G\backslash\cg_E$ then $E/_{[\alpha ]}G$ will denote the 
isomorphism class of the corresponding collapsed bundle and 
$\cg_{E/_{[\alpha ]}G}$ its gauge group.

2. The representation of $G$ on $L^2(P;E)$ generalizes the single 
scalar-valued regular representation on $L^2(G;\C )$ in two directions: it is 
a direct integral of a {\it family} of {\it vector-valued} regular 
representations. For 
instance assume that the vector bundle is globally trivial over $P$ and pick 
a trivialization $\tau :E\rightarrow P\times\C^k$. It induces an action 
$\gamma_\tau :=(\gamma\: ,\:{\rm Id}_{\C^{{\rm rk}E}})\in
C^\infty (M;{\rm Aut}P)\times{\rm Aut}\C^{{\rm rk}E}$ of $G$ on $E$ and
all actions are of this form. One can see that with respect to this
``trivial $G$-action'' (\ref{projekcio}) just reduces to ${\rm rk}_\C E$ 
copies of fiberwise integrals like (\ref{egyutthatok}) hence 
(\ref{globalis.fi.fourier}) gives back ${\rm rk}_\C E$-times the expansion 
(\ref{fourier}) parameterized by the manifold $M$.

3. For clarity we remark that the Fourier modes
(\ref{egyutthatok}) or (\ref{projekcio}) are of course independent of  
${\rm Vol}\:G$ (hence this volume will be taken from now on to be unit for 
instance).
\end{remark}


\section{The two dimensional sphere}
\label{five}


After these preliminaries providing the general framework, let us
focus our attention first to the case of the two-sphere. It will serve as 
a trivial warming-up exercise. We will demonstrate that $S^2$ admits both 
integrable strong and weak spontaneous symmetry breakings 
${\rm SO}(2)\Rrightarrow{\rm U}(1)$ and ${\rm SO}(2)\Rightarrow{\rm U}(1)$ 
respectively and of course these coincide by uniqueness hence both yield 
$S^2\cong\C P^1$. This follows from the special 
isomoprhism ${\rm SO}(2)\cong{\rm U}(1)$. 
 
First consider an {\it integrable strong spontaneous symmetry breaking} 
${\rm SO}(2)\Rrightarrow{\rm U}(1)$ over $S^2$. Put the standard round 
metric $g$ onto $S^2$ and let $\nabla$ be the corresponding Levi--Civita 
connection. Fixing an orientation on $S^2$ as well we obtain a Hodge 
operator $*:\wedge^1S^2\rightarrow\wedge^1S^2$ satisfying 
$*^2=-{\rm Id}_{\wedge^1S^2}$. By the aid of the metric 
$\wedge^1S^2\cong TS^2$ therefore we come up with an almost complex tensor 
$I:TS^2\rightarrow TS^2$. Since only the metric was used to construct it, we 
know that $\nabla I=0$ moreover it is orthogonal for $g$. Putting these 
data $(\nabla, I,g)$ into the corresponding Lagrangian we find
\[\ce (\nabla , I, g)=\frac{1}{2}\int\limits_{S^2}\left(
\frac{1}{e^2}\vert{\mathfrak S}(R_{\nabla})\vert^2_g+
\vert\nabla g\vert^2_g+\vert \nabla I\vert^2_g+
\vert II^*-{\rm Id}_{TS^2}\vert^2_g\right)\dd V =0\]
in other words $(\nabla , I,g)$ gives rise to vacuum solution to the 
strong integrable spontaneous symmetry breaking 
${\rm SO}(2)\Rrightarrow{\rm U}(1)$ yielding the usual K\"ahler structure 
on $S^2$ identifying it with $\C P^1$.

Secondly consider an {\it integrable weak spontaneous symmetry 
breaking} ${\rm SO}(2)\Rightarrow{\rm U}(1)$ over $S^2$. Because here we 
want to motivate the construction of Sect. \ref{four} step-by-step we should 
begin with the quotient ${\rm SO}(3)/{\rm U}(1)\cong S^2$ given by the 
identifications ${\rm SO}(3)\cong{\rm Aut}\HH$ and $S^2\subset{\rm Im}\HH$. 
However ${\rm SO}(3)$ is not an even dimensional compact Lie group hence it 
does not admit a complex structure {\it \`a la} Samelson therefore we even 
cannot make the very first step. To avoid this we will rather use another 
quotient ${\rm SO}(4)/{\rm U}(2)\cong S^2$. Since the following construction 
will just look like an overcomplicated version of the previous one we will 
only sketch its main steps and refer to Sect. \ref{four} for the full details.

Consider $\R^4$ equipped with an orientation and scalar product providing us 
with the group ${\rm SO}(4)$. Taking any identification $\R^4\cong\C^2$ we 
obtain a subgroup ${\rm U}(2)\subset{\rm SO}(4)$. Then ${\rm SO}(4)$ acts 
on itself transitively from the left as well as ${\rm U}(2)\subset{\rm SO}(4)$ 
acts from the right. Dividing by this right action 
we obtain a principal ${\rm U}(2)$-bundle $\pi :{\rm SO}(4)\rightarrow S^2$. 
A root decomposition of ${\mathfrak s}{\mathfrak o}(4)$ with respect to one 
of its special maximal torus $T\subset{\rm U}(2)\subset{\rm SO}(4)$ gives 
rise to a well-defined left-${\rm SO}(4)$-invariant splitting 
\[T{\rm SO}(4)=V\oplus H\]
such that $V\vert_{\pi^{-1}(x)}=T\pi^{-1}(x)\cong T{\rm U}(2)$ for all 
$x\in S^2$ moreover $H$ is also fixed by the root decomposition. The 
complexified adjoint representation of ${\rm U}(2)$ on 
${\mathfrak s}{\mathfrak o}(4)^\C$ is reducible and decomposes as 
\[{\mathfrak s}{\mathfrak o}(4)^\C\cong{\mathfrak u}(2)^\C
\oplus\C\oplus\overline\C\] 
where ${\mathfrak u}(2)^\C$ is the complexified adjoint representation of 
${\rm U}(2)$ (still reducible!) and $\C$ is the standard $1$ 
dimensional complex representation of a certain 
${\rm U}(1)\subset{\rm U}(2)$ and $\overline\C$ is its complex conjugate 
representation respectively. This gives rise to a refined left-invariant 
splitting 
\[(T{\rm SO}(4))^\C=V^\C\oplus L\oplus\overline{L}\]
such that the smooth complex line bundle $L$ over ${\rm SO}(4)$ satisfies 
$H^\C\cong L\oplus\overline{L}$ where $H^\C$ is the complexified 
horizontal bundle of the original splitting. We put two almost complex 
structures $I_H$ and $J_H$ onto $H\subset T{\rm SO}(4)$ as follows. 
The first one, $I_H: H\rightarrow H$ is defined by requiring its 
complex linear extension $I^\C_H: H^\C\rightarrow H^\C$ to satisfy  
\[ (H,I_H)\cong H^{1,0}_I:=L\subset H^\C\]
where $H^{1,0}_I\subset H^\C$ is the $+\ii$-eigenbundle of $I^\C_H$. 
This almost complex structure is left-${\rm SO}(4)$-invariant and 
right-${\rm U}(2)$-invariant by its construction. Concerning the second one 
$J_H:H\rightarrow H$, let $J$ be an integrable almost complex structure on 
${\rm SO}(4)$ found by Samelson as in Sect. \ref{two}. Consequently it has 
vanishing Nijenhuis tensor: $N_J=0$. Moreover it turns out that it is 
blockdiagonal with respect to $T{\rm SO}(4)=V\oplus H$ hence restricts to an 
almost complex structure $J_H:H\rightarrow H$. It is by construction 
left-${\rm SO}(4)$-invariant.

Putting a left-${\rm SO}(4)$-invariant and right-${\rm U}(2)$-invariant metric 
$g_0$ onto ${\rm SO}(4)$ and taking its restriction $g_H:=g_0\vert_H$ as 
well restricting the orientation of ${\rm SO}(4)$ induced by $J$ to 
$H\subset T{\rm SO}(4)$ we obtain a gauge group $\cg_H$ consisting of 
${\rm SO}(2)$ gauge transformations of $H$. Let 
\[\nabla_H: C^\infty ({\rm SO}(4);H)\times C^\infty ({\rm SO}(4);H)
\stackrel{\nabla_0\vert_H}{\longrightarrow} 
C^\infty ({\rm SO}(4);T{\rm SO}(4))\stackrel{P_H}{\longrightarrow} 
C^\infty ({\rm SO}(4);H)\] 
be the restricted-projected connection where $P_H$ is the 
$g_0$-orthogonal bundle projection from $T{\rm SO}(4)$ onto $H$. Since 
both $I_H$ and $J_H$ are left-${\rm SO}(4)$-invariant there exists a gauge 
transformation $\alpha\in\cg_H$ which rotates $I_H$ into $J_H$ i.e., 
$J_H =\alpha I_H\alpha^{-1}$. But if ${\rm U}(1)\subseteqq{\rm SO}(2)$ 
denotes the group induced by $I_H$ then taking into account that in fact 
${\rm SO}(2)\cong{\rm U}(1)$ we find that $J_H = I_H$ as complex structures 
on $H$.

One can immeditately see from our construction so far 
(especially from the fact that $J_H=I_H$)---or can Fourier 
expand $J_H$ and $g_H$ with respect to 
$[{\rm Id}_H]\in\cg_{{\rm U}(2)}\backslash\cg_H$ (cf. Definition 
\ref{ekv.fourier})---that the constructed data $(\nabla_H, J_H, g_H)$ descend 
to $S^2$ more precisely to 
\[H/_{[{\rm Id}_H]}{\rm U}(2)=TS^2.\]
The resulting triple $(\nabla ,J,g)$ on $TS^2$ is nothing else than the 
standard almost complex structure $I$ and the round metric $g$ with its 
Levi--Civita connection $\nabla$ hence it satisfies 
\[\ce (\nabla ,J ,g)=\frac{1}{2}\int\limits_{S^2}\left(
\frac{1}{e^2}\vert {\mathfrak S}(R_\nabla )\vert^2_g+
\vert\nabla g\vert^2_g+\vert N_\nabla J\vert^2_g+e^2\vert JJ^*-
{\rm Id}_{TS^2}\vert^2_g\right)\dd V=0\]
providing us with an integrable weak spontaneous symmetry breaking ${\rm SO}(2)
\Rightarrow{\rm U}(1)$. Of course because $J=I$ on $S^2$ we conclude that the 
resulting complex structure just coincides with the standard one i.e., 
$S^2\cong\C P^1$ again. 

For the technical details we refer to a completely analogous (but 
technically more complicated) construction of Sect. \ref{four}.


\section{The six dimensional sphere}
\label{four}


The time has come to carefully perform an integrable weak 
spontaneous symmetry breaking procedure over the six-sphere.

To begin with, we ask if an {\it integrable strong spontaneous symmetry 
breaking} ${\rm SO}(6)\Rrightarrow{\rm U}(3)$ exists over $S^6$. First let us 
recall the topological origin why almost complex structures exist on the 
six-sphere. Take the standard representation of SO$(6)$ on $\R^6$ and let $E$ 
be an SO$(6)$ vector bundle on $S^6$ associated to this representation. Up to 
isomorphism, these bundles are classified by the homotopy equivalence 
classes of maps from the equator of $S^6$ into SO$(6)$ i.e., by 
$\pi_5({\rm SO}(6))\cong\Z$. In particular picking an isomorphism 
$E\cong TS^6$ is equivalent to the existence of an orientation and a 
Riemannian metric on $S^6$. However a continuous map 
$f: S^5\rightarrow {\rm SO}(6)$ is homotopic to a 
continuous map $f': S^5\rightarrow {\rm U}(3)\subset{\rm SO}(6)$ because any 
particular embedding $i:{\rm U}(3)\subset {\rm SO}(6)$ induces a 
homomorphism $i_*:\pi_5({\rm U}(3))\rightarrow\pi_5({\rm SO}(6))$
which is an isomorphism: since ${\rm SO}(6)/{\rm U}(3)\cong\C P^3$ and 
$\pi_k(\C P^3)\cong 0$ if $k=5,6$ the desired isomorphism is provided by the 
associated homotopy exact sequence of the ${\rm U}(3)$-fibration 
${\rm SO}(6)\rightarrow\C P^3$. Therefore in fact any ${\rm SO}(6)$ vector 
bundle isomorphism $E\cong TS^6$ descends in a non-unique way to an U$(3)$ 
vector bundle isomorphism $E'\cong TS^6$ which is equivalent to saying that 
in addition to an orientation and a Riemannian metric $g$ on $S^6$ there 
exists a further (non-unique) compatible almost complex structure 
$J:TS^6\rightarrow TS^6$. That is $J^2=-{\rm Id}_{TS^6}$ in a manner such 
that it is orthogonal with respect to $g$ which means that 
$g(Ju, Jv)=g(u,v)$ for all vector fields $u,v$ on $S^6$. One can assume 
that these structures are smooth.

A fixed $g$ and $J$ uniquely determine a smooth non-degenerate 
$2$-form by $\omega (u,v):=g(Ju, v)$ which is not closed in general. 
In fact it cannot be closed at all: 
$\int_{S^6}\omega\wedge\omega\wedge\omega\sim{\rm Vol}(S^6)>0$ on 
the one hand but since $H^2(S^6;\Z )=0$ if $\dd\omega =0$ then $\omega 
=\dd\xi$ for a $1$-form $\xi$ hence one would find 
$\int_{S^6}\omega\wedge\omega\wedge\omega 
=\int_{S^6}\dd(\xi\wedge\dd\xi\wedge\dd\xi )=0$ via Stokes' theorem 
on the other hand, a contradiction. Consequently with any choice of $g$ and 
$J$ we cannot expect the metric to be K\"ahler or even almost K\"ahler. The 
quadruple $(S^6, g,J,\omega )$ is called an {\it almost Hermitian structure} 
on $S^6$. 

The non-existence of a K\"ahler structure on $S^6$ can be 
reformulated as saying that for any $(S^6, g,J,\omega )$ the almost 
complex tensor is not parallel with respect to the connection $\nabla$ on the 
endomorphism bundle, associated to the Levi--Civita connection i.e., 
$\nabla J\not=0$. Or using our physical language: {\it no integrable 
strong spontaneous symmetry breaking ${\rm SO}(6)\Rrightarrow{\rm U}(3)$ 
exists over the six-sphere.} 

Therefore we proceed forward ant try to construct an {\it integrable weak 
spontaneous symmetry breaking} ${\rm SO}(6)\Rightarrow{\rm U}(3)$ over $S^6$. 
Consider a Yang--Mills--Higgs--Nijenhuis field theory (\ref{YMH}) specialized 
to the six-sphere:
\begin{equation}
\ce (\nabla ,\Phi ,g):=\frac{1}{2}\int\limits_{S^6}\left(
\frac{1}{e^2}\vert {\mathfrak S}(R_\nabla )\vert^2_g+
\vert\nabla g\vert^2_g+
\vert N_\nabla\Phi\vert^2_g+e^2\vert\Phi\Phi^*-
{\rm Id}_{TS^6}\vert^2_g\right)\dd V. 
\label{S6YMH}
\end{equation}
Recall that the existence of a smooth solution to the corresponding vacuum 
equations (\ref{vakuum}) is equivalent to the existence of a non-K\"ahlerian 
complex structure on $S^6$. Our aim is therefore to prove the existence of a 
spontaneously broken vacuum in this theory. This task will be carried out by 
the aid of the Lie group ${\rm G}_2$ and its complex structure in several 
steps below. 

So let us consider the $14$ dimensional connected, compact, simply connected 
exceptional Lie group ${\rm G}_2$. It acts on itself by 
{\it left}-translations. As it is well-known in addition letting ${\rm SU}(3)$ 
act on ${\rm G}_2$ from the {\it right} ${\rm G}_2$ arises as the total space 
of an ${\rm SU}(3)$ fibration $\pi : {\rm G}_2\rightarrow S^6$. More precisely 
(for a complete proof cf. e.g. \cite{bry},\cite[pp. 306-311]{pos}) taking the 
identifications ${\rm G}_2\cong{\rm Aut}\OO$ and $S^6\subset{\rm Im}\OO$ the 
group ${\rm G}_2$ acts on $S^6$ transitively and the corresponding isotropy 
subgroup can be found as follows. Over a point $x\in S^6\subset{\rm Im}\OO$ 
of unit length---with respect to the 
standard scalar product---the tangent space $T_xS^6$ can be identified with a 
subspace $V_x\subset\OO$ perpendicular to the $\R$-span of $\{ 1, x\}$ in 
$\OO$. Orientation and the scalar product on $\OO$ restricts to $V_x$ as well 
as a complex structure on $V_x$ is induced by Cayley multiplication with $x$ 
itself rendering $V_x$ an oriented $3$ dimensional complex vector space 
with a Hermitian scalar product. A fiber of our principal bundle over 
$x\in S^6$ then arises as the special unitary group of this $V_x$. This 
principal bundle is non-trivial: the group ${\rm G}_2$ as a principal 
${\rm SU}(3)$-bundle over $S^6$ corresponds to a generator of 
$\pi_5({\rm SU}(3))\cong\Z$, cf. e.g. \cite[Corollary 3]{cha-rig2}. 

First referring to this fibration and the Lie algebra structure on the 
tangent bundle of ${\rm G}_2$ we construct a left-invariant splitting  
\[T{\rm G}_2=V\oplus H\] 
as follows. For the vertical bundle put 
$V\vert_{\pi^{-1}(x)}:=T\pi^{-1}(x)\cong T{\rm SU}(3)$ for all $x\in S^6$.
Regarding the horizontal bundle $H\cong (T{\rm G}_2)/V$
note that a positive root system $P^+$ of ${\rm SU}(3)$
can be a subsystem of a positive root system $R^+$ of ${\rm G}_2$ hence 
with the special Cartan subalgebra ${\mathfrak t}\subset {\mathfrak s}
{\mathfrak u}(3)\subset{\mathfrak g}_2$ the decomposition of
${\mathfrak g}_2$ into real root spaces (\ref{gyokok}) yields ${\mathfrak g}_2=
{\mathfrak s}{\mathfrak u}(3) \oplus{\mathfrak m}$ where
\[{\mathfrak s}{\mathfrak u}(3)={\mathfrak t}\oplus\bigoplus
\limits_{\mu\in P^+}{\mathfrak g}_\mu,\:\:\:\:\: 
{\mathfrak m}:=\bigoplus\limits_{\nu\in R^+\setminus P^+}
{\mathfrak g}_\nu\]
fixing the choice of ${\mathfrak m}\subset {\mathfrak g}_2$. Identifying
both ${\mathfrak s}{\mathfrak u}(3)$ and ${\mathfrak m}$ with subspaces of
left-invariant sections we fix $T{\rm G}_2= V\oplus H$ up to a Cartan 
subalgebra of ${\rm SU}(3)$ as follows: we already know that 
$V_y={\rm ev}_y{\mathfrak s}{\mathfrak u}(3)$ moreover put 
$H_y:={\rm ev}_y{\mathfrak m}$ for all $y\in {\rm G}_2$. In particular $H$ is 
a trivial bundle over ${\rm G}_2$. 

The complexification of this splitting can be refined. The complexified 
adjoint representation of ${\rm SU}(3)$ on ${\mathfrak g}_2^\C$ is 
reducible and its decomposition into irreducible summands is given by
\[{\mathfrak g}_2^\C\cong {\mathfrak s}{\mathfrak u}(3)^\C\oplus
\C^3\oplus\overline{\C^3}\]
where ${\mathfrak s} {\mathfrak u}(3)^\C$ is the complexified
adjoint representation while $\C^3$ and $\overline{\C^3}$ are the
standard $3$ dimensional complex representation of ${\rm SU}(3)$
and its complex conjugate respectively. This induces a 
refined left-invariant decomposition 
\[(T{\rm G}_2)^\C =V^\C\oplus Z\oplus\overline{Z}\]
of $(T{\rm G}_2)^\C =V^\C\oplus H^\C$ such that the smooth complex vector 
bundle $Z$ over ${\rm G}_2$ satisfies $H^\C =Z\oplus\overline{Z}$ 
where $H^\C :=H\otimes\C$ is the complexified horizontal bundle.

Secondly we consider two canonical almost complex structures $I_H$ and $J_H$ 
on $H\subset T{\rm G}_2$. The first one, $I_H: H\rightarrow H$ is defined by 
requiring that its complex linear extension $I^\C_H :H^\C\rightarrow H^\C$ 
satisfies 
\[(H, I_H)\cong H^{1,0}_I:=Z\subset H^\C\]
where $H^{1,0}_I\subset H^\C$ is the $+\ii$-eigenbundle of $I_H^\C$. Note 
that by construction $I_H: H\rightarrow H$ is a left-${\rm G}_2$-invariant 
and right-${\rm SU}(3)$-invariant almost complex structure which is 
moreover not integrable.\footnote{We note that by right 
${\rm SU}(3)$-invariance it descends and gives the Cayley almost complex 
structure on $S^6$.}

The other one, $J_H:H\rightarrow H$ is constructed as follows. We already 
know from Sect. \ref{two} that there exists a left-${\rm G}_2$-invariant 
integrable almost complex tensor $J$ on ${\rm G}_2$ {\it \`a la} Samelson. 
It has of course vanishing Nijenhuis tensor 
\begin{equation}
N_J=0
\label{zero}
\end{equation}
moreover via (\ref{J-szetszedes}) it admits a left-invariant
global blockdiagonal decomposition
\begin{equation}
J=\begin{pmatrix}
                  J_V& 0\\
                   0 & J_H
         \end{pmatrix}
\label{blokk}
\end{equation}  
where $J_V\in C^\infty ({\rm G}_2;{\rm End}V)$ and $J_H\in 
C^\infty ({\rm G}_2;{\rm End}H)$. In particular $J_H^2=-{\rm Id}_H$ and in 
this way we obtain a smooth trivial rank $3$ complex vector bundle $(H, J_H)$ 
over ${\rm G}_2$. That is, as a complex vector bundle 
\[(H, J_H)\cong H^{1,0}_J\subset H^\C\] 
where $H^{1,0}_J$ is the $+\ii$-eigenbundle of $J_H^\C: H^\C\rightarrow H^\C$. 

Thirdly we construct an invariant metric on ${\rm G}_2$. Fix an 
orientation on ${\rm G}_2$ induced by $J$. Also fix a Riemannian metric 
$g_0$ on ${\rm G}_2$ which is invariant under both 
left-${\rm G}_2$-translations induced by the left-action of ${\rm G}_2$ and
fiberwise right-${\rm SU}(3)$-translations induced by the right-action of
${\rm SU}(3)\subset{\rm G}_2$. Such metric $g_0$ easily arises by averaging an 
arbitrary metric $g$ on ${\rm G}_2$ with respect to ${\rm G}_2\times
{\rm SU}(3)$ equipped with the measure $\dd y\otimes\dd\gamma$
where $\dd y$ and $\dd\gamma$ are the unit-volume bi-invariant Haar measures 
on ${\rm G}_2$ and ${\rm SU}(3)$ respectively:
\[g_0(u,v):=\int\limits_{{\rm G}_2\times{\rm SU}(3)}
g\left( (L_y)_*(R_\gamma )_*u\:,\:(L_y)_*(R_\gamma )_*v\right)
\dd y\otimes\dd\gamma .\]
If the metric is normalized such that the induced volume form satisfies 
$\int_{{\rm G}_2}\dd V_0=1$ then as left-invariant measures $\dd V_0=\dd y$. 
The associated Levi--Civita connection will be denoted by $\nabla_0$. 
Taking the restricted orientation and metric $g_0\vert_H=:g_H$ we also 
have an associated gauge group $\cg_H$ consisting of ${\rm SO}(6)$ gauge 
transformations of the horizontal bundle. Take the restricted-projected 
connection 
\begin{equation}
\nabla_H: C^\infty ({\rm G}_2;H)\times C^\infty ({\rm 
G}_2;H)\stackrel{\nabla_0\vert_H}{\longrightarrow}C^\infty ({\rm 
G}_2;T{\rm G}_2)\stackrel{P_H}{\longrightarrow}C^\infty ({\rm G}_2;H)
\label{megszoritott.konnexio}
\end{equation} 
where $P_H$ is the left-${\rm G}_2$-equivariant $g_0$-orthogonal bundle 
projection from $T{\rm G}_2$ onto $H$. The metric furthermore admits a unique 
Hermitian extension $g_H^\C$ to $H^\C\cong H\oplus\ii H$ given by
\[g_H^\C (u_1+\ii u_2, v_1+\ii v_2):=
g_H(u_1,v_1)+g_H(u_2,v_2)+\ii\left( g_H(u_2,v_1)-g_H(u_1,v_2)\right)\]
(complex conjugate linear in its second variable) yielding the complexified 
gauge group $\cg^\C_H$ consisting of ${\rm U}(6)$ gauge transformations with 
respect to $g_H^\C$.

By the aid of this gauge group the almost complex structures $J_H$ and 
$I_H$ on the horizontal bundle can be compared. By construction 
$I_H$ commutes with the fiberwise right-action of ${\rm SU}(3)
\subset{\rm G}_2$ hence for all $y\in{\rm G}_2$ and $\gamma\in{\rm SU}(3)$ 
we know that
\begin{equation}
I_{H_{y\gamma}}(R_\gamma )_*=(R_\gamma )_*I_{H_y}.
\label{folcsereles}
\end{equation}
In addition $I_H$ is invariant under the left-${\rm G}_2$-action hence 
$H^{1,0}_I\subset H^\C$ is a left-invariant bundle moreover $I_H$ is 
orthogonal for $g_H$. But $J_H$ is also invariant under the 
left-${\rm G}_2$-action hence $H^{1,0}_J\subset H^\C$ is also left-invariant by 
construction. Therefore these two sub-bundles can be rotated into 
each other within $H^\C$ i.e., there exists a gauge transformation 
$\alpha\in\cg^\C_H$ such that $H^{1,0}_J=\alpha H^{1,0}_I$ i.e., there 
exists an induced real gauge transformation $\alpha\in\cg_H$ such 
that for all $y\in {\rm G}_2$
\[J_{H_y}=\alpha (y)I_{H_y}\alpha^{-1}(y)\]
holds. Note that in any point 
$\alpha (y)\notin{\rm U}(3)$ with ${\rm U}(3)\subset{\rm SO}(6)$ being induced 
by $I_{H_y}$ because $I_{H_y}$ and $J_{H_y}$ are not equivalent; but at least 
$\alpha (y)\in{\rm SO}(6)$ where ${\rm SO}(6)$ is by definition induced by 
$g_H$ and the orientation.\footnote{In our terminology developed at the end of
Sect. \ref{two} the element $\alpha\in\cg_H$ acts as a {\it deformation} of
$J_H$ into $I_H$.} Hence in particular it follows that $J_H$ is also 
orthogonal with respect to $g_H$. 

It follows from (\ref{blokk}) that $J_HJ_H^*={\rm Id}_H$ and 
(\ref{lokalis-nijenhuis}) together with (\ref{blokk}) also imply that 
$N_{\nabla_H} J_H=N_{\nabla_0\vert_H}J_H$. Moreover we deduce from
(\ref{megszoritott.konnexio}) that $\vert\nabla_Hg_H\vert_{g_H}\leqq\vert 
(\nabla_0\vert_H)g_H\vert_{g_H}$ and taking into account that
$\nabla_0\vert_H$ is torsion-free we also find ${\mathfrak S}(R_{\nabla_H})=
{\mathfrak S}(P_HR_{\nabla_0\vert_H})$ hence $\vert 
{\mathfrak S}(R_{\nabla_H})\vert_{g_H}\leqq
\vert {\mathfrak S}(R_{\nabla_0\vert_H})\vert_{g_H}$. Therefore 
\begin{eqnarray}
\ce (\nabla_H, J_H,g_H)&=&\frac{1}{2}\int\limits_{{\rm G}_2}
\left(\frac{1}{e^2}\vert {\mathfrak S}(R_{\nabla_H})\vert^2_{g_H}+
\vert\nabla_Hg_H\vert^2_{g_H}+\vert N_{\nabla_H} J_H\vert^2_{g_H}+e^2\vert
J_HJ_H^*-{\rm Id}_H\vert^2_{g_H}\right)\dd V_0\nonumber\\
&\leqq&\mbox{ horizontal part of $\ce (\nabla_0, J, g_0)$}
\leqq\ce (\nabla_0, J, g_0)\nonumber
\end{eqnarray}
and via (\ref{zero}) the functional (\ref{YMH}) specialized to ${\rm G}_2$
satisfies $\ce (\nabla_0, J, g_0)=0$. Our constructions so far can be 
then summarized as follows: the triple $(\nabla_H, J_H, g_H)$ represents a 
smooth vacuum solution (\ref{vakuum}) to the Yang--Mills--Higgs--Nijenhuis 
field theory (\ref{YMH}) restricted to the horizontal bundle 
$H\subset T{\rm G}_2$ i.e., 
\begin{equation}
\ce(\nabla_H, J_H, g_H)=0.
\label{G2vakuum}
\end{equation}
We want to use this vacuum solution to obtain a vacuum solution to 
(\ref{S6YMH}). Our technical tool to achieve this will be a carefully 
constructed global Fourier expansion with respect to ${\rm SU}(3)$ of the 
vacuum Higgs field $J_H$.

As we have seen in Sect. \ref{three} for this purpose we need an action 
of ${\rm SU}(3)$ on the Hilbert space of sections of the horizontal bundle 
$H\subset T{\rm G}_2$. Given any $\gamma\in{\rm SU}(3)$ and $u_y\in H_y$ 
consider the real linear isomorphism 
$\gamma_\alpha :H_y\rightarrow H_{y\gamma}$ whose shape is
\[u_y\gamma_\alpha :=\left(\pi_*\vert_{H_{y\gamma}}\alpha^{-1}(y\gamma )
\right)^{-1}\left(\pi_*\vert_{H_y}\alpha^{-1}(y)\right) u_y\]
or in other words, making use of the fiberwise right ${\rm SU}(3)$-translation
$R_\gamma :{\rm G}_2 \rightarrow {\rm G}_2$ we put
\begin{equation}
u_y\gamma_\alpha :=
\left(\alpha (y\gamma )(R_\gamma )_*\alpha^{-1}(y)\right)u_y.
\label{B}
\end{equation}
Dividing by this $\alpha$-dependent action we get an abstract isomorphism 
$H/_{[\alpha ]}{\rm SU}(3)\cong TS^6$ for the quotient as a plain real 
vector bundle.

The space $C^\infty ({\rm G}_2; H)$ carries an associated representation 
of ${\rm SU}(3)$ from the left defined by
\begin{equation}
(\gamma_\alpha\cdot u)(y):=u(y\gamma )\gamma^{-1}_\alpha
\label{abrazolas}
\end{equation}
for all $u\in C^\infty ({\rm G}_2;H)$. Consider the real Hilbert space 
$L^2({\rm G}_2;H)$ as the completion of $C^\infty ({\rm G}_2; H)$ 
with respect to the $L^2$ scalar product $(u,v)_{L^2({\rm G}_2)}:=
\int_{{\rm G}_2}g_H(u(y),v(y))\dd V_0$. 
Extending (\ref{abrazolas}) from $C^\infty ({\rm G}_2;H)$ to this space we 
obtain continuous real representations of ${\rm SU}(3)$ parameterized by 
$\cg_H$. These representations are orthogonal because the metric $g_H$ 
is right-invariant with respect to ${\rm SU}(3)$. The construction can be 
complexified and we come up with continuous complex unitary 
representations of ${\rm SU}(3)$ on the fixed complex Hilbert space 
$L^2({\rm G}_2;H^\C )\cong L^2({\rm G}_2;H)\otimes\C$ parameterized by the 
complexified gauge group $\cg^\C_H$. All of these representations are 
unitarily equivalent.

Now we apply the general construction of Sect. \ref{three} to our 
situation by taking the compact group $G$ to be ${\rm SU}(3)$, the 
total space $P$ to be ${\rm G}_2$ and the auxiliary vector bundle $E$ on it to 
be ${\rm End}H$. Note that $J_H\in C^\infty 
({\rm G}_2\:;\:{\rm End}H)$ because 
${\rm Ad}{\rm G}_2\subset {\rm End}\:T{\rm G}_2$. Picking a gauge 
transformation $\alpha\in\cg_H$ there is an action of ${\rm SU}(3)$ on 
$L^2({\rm G}_2;{\rm End}H)$ induced by (\ref{abrazolas}). 
Consequently the smooth Higgs field $J_H$ admits an associated 
$\cg_{{\rm SU}(3)}$-equivariant global Fourier expansion 
(\ref{globalis.fi.fourier}) of the shape
\begin{equation}
J_H =\sum\limits_{\rho\in{\rm Irr}({\rm SU}(3);\C )}\Pi_{\alpha ,\rho} J_H
\label{sorfejtes}
\end{equation}
and here pointwise equality holds.

\begin{lemma}
Let $J_H$ be the horizontal part of a left-${\rm G}_2$-invariant 
integrable and let $I_H$ be the horizontal part of a 
left-${\rm G}_2$-invariant and right ${\rm SU}(3)$-invariant non-integrable 
almost complex structure on ${\rm G}_2$ respectively. Consider a gauge 
transformation $\alpha\in\cg_H$ with respect to the metric $g_H$ 
satisfying $J_H=\alpha I_H\alpha^{-1}$ as above. 

Then the global Fourier expansion (\ref{sorfejtes}) corresponding 
to $\alpha\in\cg_H$ satisfies $J_H =\Pi_{\alpha ,1} J_H$ 
i.e., all the higher Fourier modes vanish: $\Pi_{\alpha ,\rho} J_H =0$ with 
$\rho\not\cong 1$. Consequently $J_H$ uniquely descends to an endomorphism 
\[J_{H/_{[\alpha ]}{\rm SU}(3)}: H/_{[\alpha ]}{\rm SU}(3)\longrightarrow 
H/_{[\alpha ]}{\rm SU}(3)\] 
over ${\rm G}_2/{\rm SU}(3)\cong S^6$.

Moreover this Fourier expansion is well-defined up to an element of 
$\cg_{{\rm SU}(3)}\subset\cg_H$ (cf. Definition \ref{ekv.fourier}) i.e., 
depends only on the class $[\alpha ]\in \cg_{{\rm SU}(3)}\backslash\cg_H$. 
\label{vakuumlemma}
\end{lemma}
\noindent {\it Proof.} Exploiting the identities (\ref{folcsereles}) and 
(\ref{B}) we calculate the ground mode $\Pi_{\alpha ,1}J_H$ in 
(\ref{sorfejtes}) over a point $y\in{\rm G}_2$ as follows:
\begin{eqnarray}
\Pi_{\alpha ,1}J_{H_y}& =& \int\limits_{{\rm SU}(3)}\left(\gamma_\alpha\cdot
J_{H_y}\right)\dd\gamma\nonumber\\
         &=& \int\limits_{{\rm SU}(3)}\left(\gamma^{-1}_\alpha
J_{H_{y\gamma}}\gamma_\alpha\right)\dd\gamma\nonumber\\
        &=&\int\limits_{{\rm SU}(3)}\left(\gamma^{-1}_\alpha \left(\alpha
(y\gamma)I_{H_{y\gamma}}\alpha^{-1}(y\gamma)\right)\gamma_\alpha\right)
\dd\gamma\nonumber\\
& = &\int\limits_{{\rm SU}(3)}\left(\alpha (y)(R_{\gamma^{-1}})_*
\alpha^{-1}(y\gamma)\alpha (y\gamma)I_{H_{y\gamma}}\alpha^{-1}(y\gamma)\alpha
(y\gamma )(R_\gamma )_*\alpha^{-1}(y)\right)\dd\gamma\nonumber\\
& =& \int\limits_{{\rm SU}(3)}\left(\alpha (y)
(R_{\gamma^{-1}})_*I_{H_{y\gamma}}(R_\gamma)_*\alpha^{-1}(y)\right)
\dd\gamma\nonumber\\
&=&\int\limits_{{\rm SU}(3)}\left(\alpha (y)I_{H_y}\alpha^{-1}(y)\right)
\dd\gamma\nonumber\\
&=&\int\limits_{{\rm SU}(3)}J_{H_y}\dd\gamma\nonumber\\
&=&J_{H_y}.\nonumber
\end{eqnarray}
Therefore completeness gives $\Pi_{\alpha ,\rho} J_H=0$ if $\rho\not\cong 1$ 
and we obtain the result. 

Moreover this result obviously remains unchanged if $\alpha\in\cg_H$ is 
multiplied from the left with a vertically constant gauge transformation 
$\beta\in\cg_{{\rm SU}(3)}\subset\cg_H$ as claimed. $\Diamond$
\vspace{0.1in}

\begin{lemma} Making use of the gauge transformation $\alpha\in\cg_H$ 
of Lemma \ref{vakuumlemma} and the associated Fourier expansion the 
horizontal part $g_H$ of the metric on $H\subset T{\rm G}_2$ satisfies 
$g_H=\Pi_{\alpha, 1}g_H$. Consequently $g_H$ uniquely descends to a 
metric $g_{H/_{[\alpha ]}{\rm SU}(3)}$ on the bundle $H/_{[\alpha ]}{\rm 
SU}(3)$ over ${\rm G}_2/{\rm SU}(3)\cong S^6$. 

Additionally, the associated connection $\nabla_H$ in 
(\ref{megszoritott.konnexio}) uniquely descends to a connection
\[\nabla_{H/_{[\alpha ]}{\rm SU}(3)} :
C^\infty (S^6; H/_{[\alpha ]}{\rm SU}(3))\longrightarrow
C^\infty (S^6; (H/_{[\alpha ]}{\rm SU}(3))\otimes\wedge^1S^6)\]
compatible with $g_{H/_{[\alpha ]}{\rm SU}(3)}$ and gives rise to the
corresponding curvature tensor
\[R_{\nabla_{H/_{[\alpha ]}{\rm SU}(3)}} :C^\infty (S^6;
H/_{[\alpha ]}{\rm SU}(3))\longrightarrow
C^\infty (S^6; (H/_{[\alpha ]}{\rm SU}(3))\otimes\wedge^2S^6).\]
\label{vakuumlemma'}
\end{lemma}

\noindent{\it Proof.} The metric $g_H$ is 
right-${\rm SU}(3)$-translation invariant: for any $u,v\in
C^\infty ({\rm G}_2;H)$ one finds $g_H(u,v)=
g_H((R_\gamma )_*u, (R_\gamma )_*v)$ with $\gamma\in{\rm SU}(3)$; as 
well as it is of course invariant under its own gauge group: 
$g_H(u,v)=g_H(\alpha u, \alpha v)$ with $\alpha\in\cg_H$. Therefore by 
(\ref{B}) we find that
\begin{equation}
g_H(u,v)=g_H(u\gamma_\alpha ,v\gamma_\alpha )
\label{invariancia}
\end{equation} 
i.e., it is invariant under any twisted right-action. In particular we can 
Fourier expand the metric with respect to the action of ${\rm SU}(3)$ on 
$L^2({\rm G}_2;H^*\otimes H^*)$ induced by (\ref{abrazolas}) with 
$\alpha\in\cg_H$ used in Lemma \ref{vakuumlemma}. Taking account 
(\ref{invariancia}) a similar calculation as in Lemma \ref{vakuumlemma} 
demonstrates that 
\[g_H=\Pi_{\alpha ,1}g_H\] 
hence $g_H$ descends to a metric $g_{H/_{[\alpha ]}{\rm SU}(3)}$ on the 
bundle $H/_{[\alpha ]}{\rm SU}(3)$ over $S^6$. 

Moreover the compatibility equation $\nabla_Hg_H=0$ says for arbitrary 
horizontal vector fields $u,v,w\in C^\infty ({\rm G}_2; H)$ that 
\begin{equation}
\dd\left( g_H(u,v)\right) w=g_H\left( 
(\nabla_H)_w\:u\:,\:v)+g_H(u\:,\:(\nabla_H)_w\:v\right) .
\label{kompatibilitas1}
\end{equation} 
Therefore  
\begin{eqnarray}
\dd\left( g_H(u\gamma_\alpha ,v\gamma_\alpha )\right) w
&=&
g_H\left( (\nabla_H)_w(u\gamma_\alpha )\:,\:v\gamma_\alpha\right)
+g_H\left(u\gamma_\alpha \:,\:(\nabla_H)_w(v\gamma_\alpha )\right)\nonumber\\
&&\nonumber\\
&=& 
g_H\left( (\nabla_H)_w(u\gamma_\alpha )\gamma^{-1}_\alpha\gamma_\alpha\:,\:
v\gamma_\alpha\right) +g_H\left( u\gamma_\alpha\:,\:
(\nabla_H)_w(v\gamma_\alpha )\gamma^{-1}_\alpha\gamma_\alpha\right)\nonumber\\
&&\nonumber\\
&=&
g_H\left( (\nabla_H)_w(u\gamma_\alpha )\gamma^{-1}_\alpha\:,\:v\right) +
g_H\left( u\:,\:(\nabla_H)_w(v\gamma_\alpha)\gamma^{-1}_\alpha\right) .
\label{kompatibilitas2}
\end{eqnarray}
A comparison of (\ref{invariancia}), (\ref{kompatibilitas1}) and 
(\ref{kompatibilitas2}) shows that $\nabla_H(\cdot )= 
\nabla_H((\cdot )\gamma_\alpha )\gamma^{-1}_\alpha$ i.e., $\nabla_H$ 
commutes with the twisted right-${\rm SU}(3)$-translations, too. 
Consequently we find that the restricted Levi--Civita connection 
$\nabla_H$ also descends uniquely to a connection with an induced 
curvature as claimed. $\Diamond$
\vspace{0.1in}

\noindent We proceed further and choose a particular real vector bundle 
isomorphism 
\[\varphi :TS^6\longrightarrow H/_{[\alpha ]}{\rm SU}(3).\] 
The endomorphism 
$J_{H/_{[\alpha ]}{\rm SU}(3)}: H/_{[\alpha ]}{\rm SU}(3)\rightarrow 
H/_{[\alpha ]}{\rm SU}(3)$ of Lemma \ref{vakuumlemma} can be transported to 
an almost complex structure $J:TS^6\rightarrow TS^6$ by the formula 
\[J:=\varphi^{-1}(J_{H/_{[\alpha ]}{\rm SU}(3)})\varphi .\] 
Pulling back the metric and the connection in Lemma \ref{vakuumlemma'} we 
also get a Riemannian metric
\[g:=(g_{H/_{[\alpha ]}{\rm SU}(3)})(\varphi\times\varphi )\] 
(i.e., $g(\tilde{u},\tilde{v}):=g_{H/_{[\alpha ]}{\rm SU}(3)}
(\varphi (\tilde{u}), \varphi (\tilde{v}))$ for all 
$\tilde{u},\tilde{v}\in C^\infty (S^6;TS^6)$) and a compatible connection 
\[\nabla :=\varphi^{-1}(\nabla_{H/_{[\alpha ]}{\rm SU}(3)})\varphi\] 
on $TS^6$. The triple $(\nabla, J, g)$ depends only on 
$[\alpha ]\in\cg_{{\rm SU}(3)}\backslash\cg_H$. Indeed, if 
$\psi :TS^6\rightarrow H/_{[\alpha ]}{\rm SU}(3)$ is another isomorphism then 
$\psi^{-1}\circ\varphi :TS^6\rightarrow TS^6$ is a general gauge 
transformation of the tangent bundle hence by our {\bf Principle} from Sect. 
\ref{two} we can suppose that the particular choice of the vector bundle 
isomorphism between the bundles $TS^6$ and $H/_{[\alpha ]}{\rm SU}(3)$ over 
$S^6$ does not effect the geometric structure induced by $(\nabla ,J,g)$ on 
$S^6$. In the particular case of an almost complex structure we checked
at the end of Sect. 2 that $N_\nabla J$ transforms as a tensor under the 
induced action of $\psi^{-1}\circ\varphi\in\cg_{TS^6}$. 

We make a digression here and prove the counter-intuitive fact that 

\begin{lemma} The almost complex structure $J:TS^6\rightarrow TS^6$ is not 
homogeneous.
\label{homogenlemma}
\end{lemma}

\noindent{\it Proof.} Our strategy to prove this will 
be as follows. It is well-known (cf., e.g. \cite{cal-glu}) that 
homogeneous almost complex structures on $S^6$ are parameterized by 
${\rm SO}(7)/{\rm G}_2\cong\R P^7$ and are orthogonally equivalent to 
the standard Cayley one $I$ with respect to the standard metric on $S^6$. 
Therefore if we can prove that $J$ just constructed above is not orthogonally 
equivalent to $I$ then we are done.

There exists a special isomorphism
\[\varphi_\alpha :TS^6\cong H/_{[{\rm Id}_H]}{\rm SU}(3)\longrightarrow 
H/_{[\alpha]}{\rm SU}(3)\]
induced by the gauge transformation $\alpha\in\cg_H$ from Lemma 
\ref{vakuumlemma}. For an arbitrary $\beta\in\cg_H$ let us 
denote sections of the corresponding quotient bundle as 
$[u]_\beta\in C^\infty\left( S^6; H/_{[\beta ]}{\rm SU}(3)\right)$ 
consisting of fiberwise equivalence classes of $\gamma_\beta$-invariant 
sections $u\in C^\infty ({\rm G}_2;H)$, cf. Sect. \ref{three}. Then a 
tangent vector $\tilde{u}\in C^\infty (S^6; TS^6)$ can 
be written in the form $[u]_{{\rm Id}_H}\in C^\infty
\left( S^6; H/_{[{\rm Id}_H]}{\rm SU}(3)\right)$. Since (\ref{B}) gives an 
identity $\alpha (u\gamma_{\:{\rm Id}_H})=(\alpha u)\gamma_\alpha$ we simply 
put 
\[\varphi_\alpha (\tilde{u})=\varphi_\alpha ([u]_{{\rm Id}_H})
:=[\alpha u]_\alpha\] 
with inverse $\varphi^{-1}_\alpha ([u]_\alpha )=[\alpha^{-1}u]_{{\rm Id}_H}$. 
Because elements of the form $\alpha u\in C^\infty ({\rm G}_2;H)\subset 
C^\infty ({\rm G}_2;T{\rm G}_2)$ are gauge transformed objects the almost 
complex structure $J_H:H\rightarrow H$ also acts on them in its {\it gauge
transformed} form $\alpha (J_H) = \alpha J_H\alpha^{-1}$ by our 
{\bf Principle}. Hence recalling the construction of $J$ and writing 
$J_H = \alpha I_H\alpha^{-1}$ as a usual {\it deformation} of $I_H$ (cf. 
the remarks at the end of Sect. 2) we obtain
\[J\tilde{u}=\varphi^{-1}_\alpha (J_{H/_{[\alpha]}{\rm SU}(3)}(\varphi_\alpha 
([u]_{{\rm Id}_H})))=
[\alpha^{-1}(\alpha J_H\alpha^{-1})\alpha u]_{{\rm Id}_H}=[J_Hu]_{{\rm Id}_H}=
[\alpha I_H\alpha^{-1}u]_{{\rm Id}_H}.\]
This means that the operator $\alpha I_H\alpha^{-1}$, when restricted to 
${\rm SU}(3)$-invariant sections, descends to $TS^6$ hence on this 
subspace it coincides with its ground mode in its Fourier expansion 
with respect to ${\rm Id}_H\in\cg_H$. Consequently 
$[\alpha I_H\alpha^{-1}u]_{{\rm Id}_H}$ can be calculated analytically as 
$\Pi_{{\rm Id}_H, 1}(\alpha I_H\alpha^{-1})$. Repeating the steps of the 
proof of Lemma \ref{vakuumlemma} we find for an $y\in {\rm G}_2$ that 
\[\Pi_{{\rm Id}_H, 1}\left(\alpha (y)I_{H_y}\alpha^{-1}(y)\right) =
\int\limits_{{\rm SU}(3)}\left( (R_{\gamma^{-1}})_*\alpha (y\gamma )
I_{H_{y\gamma}}\alpha^{-1}(y\gamma )(R_\gamma )_*\right)\dd\gamma .\]
Using right-{\rm SU}(3)-translation to identify 
${\rm End}(H\vert_{\pi^{-1}(x)})$ with ${\rm SU}(3)\times{\rm End}\:H_y$ 
for some fixed $y\in\pi^{-1}(x)$ we can suppose that $\alpha (y(\cdot )): 
{\rm SU}(3)\rightarrow {\rm End}\:H_y$ and $I_{H_{y(\cdot )}}: {\rm SU}(3)
\rightarrow {\rm End}\:H_y$ are functions and the latter being constant 
by its construction. Therefore we come up with
\begin{eqnarray}
\Pi_{{\rm Id}_H, 1}\left(\alpha (y)I_{H_y}\alpha^{-1}(y)\right)&=&
\int\limits_{{\rm SU}(3)}\left( \alpha (y\gamma )
I_{H_{y\gamma}}\alpha^{-1}(y\gamma )\right)\dd\gamma\nonumber\\
&=&\int\limits_{{\rm SU}(3)}\left( {\rm Ad}_{\alpha^{-1}(y\gamma )}
(I_{H_{y\gamma}})\right)\dd\gamma\nonumber\\ 
&=&\left(\:\:\int\limits_{{\rm SU}(3)}{\rm Ad}_{\alpha^{-1}(y\gamma )}
\dd\gamma\right)I_{H_y}.\nonumber
\end{eqnarray}
Introducing the real orthogonal representation 
$\sigma_y:{\rm SU}(3)\rightarrow{\rm Aut}({\rm End}H_y)$ where 
$\sigma_y(\gamma ):= {\rm Ad}_{\alpha^{-1}(y\gamma )}$ we conclude that 
\[\Pi_{{\rm Id}_H, 1}\left(\alpha (y)I_{H_y}\alpha^{-1}(y)\right)=
{\bf P}_y\left( I_{H_y}\right)\]
where ${\bf P}_y:=\int_{{\rm SU}(3)}\sigma_y(\gamma )\dd\gamma$ is a 
finite dimensional projection onto the invariant subspace of $\sigma_y$:
\[{\bf P}_y:{\rm End}\:H_y\longrightarrow ({\rm End}\:H_y)_{{\rm SU}(3)}.\]
It is easy to check that $\sigma_y$ is a reducible representation and its 
invariant subspace is spanned by two invariant operators of $H_y$ namely 
the identity operator of $H_y$ and the complex multiplication on $H_y$ 
induced by the embedding $\alpha (y{\rm SU}(3))\subset{\rm SO}(6)$. 
This complex structure on $H_y$ is obviously different from $(H_y, I_{H_y})$ 
---because $\alpha (y\gamma )\notin{\rm U}(3)\subset{\rm SO}(6)$ 
induced by $I_{H_y}$---consequently ${\bf P}_y$ indeed projects $I_{H_y}$ 
non-trivially. It is also clear from the construction that ${\bf P}_y
\left(I_{H_y}\right)$ is ${\rm SU}(3)$-invariant consequently if $\pi (y)=x$ 
then ${\bf P}_y\left(I_{H_y}\right)$ descends unambigously from $H_y$ to 
$(H/_{[{\rm Id}_H]}{\rm SU}(3))_x\cong 
T_xS^6$ and it coincides with $J_x:T_xS^6\rightarrow T_xS^6$. In other words 
there exists an element $a\in C^\infty (S^6; {\rm Aut}(TS^6))$ defined by  
\[J=:aIa^{-1}\] 
where $I$ is the Cayley almost complex structure on $S^6$. Taking into 
account that $g_H$ projects onto the standard metric $g_0$ on 
$H/_{[{\rm Id}_H]}{\rm SU}(3)\cong TS^6$ we conclude that 
the automorphism $a\in C^\infty (S^6; {\rm Aut}(TS^6))$ is a non-orthogonal 
element with respect to the standard metric $g_0$ on $S^6$. This is because 
$a (x)\in {\rm Aut}(T_xS^6)$ with $x\in S^6$ arises as the average of 
$g_H$-orthogonal transformations $\alpha (y\gamma )\in {\rm Aut}H_y$ when 
$\gamma$ runs over ${\rm SU}(3)\cong\pi^{-1}(x)\subset{\rm G}_2$ and 
orthogonality is lost during averaging. Consequently $J$ is not orthogonally 
equivalent to the Cayley structure as claimed.

Finally we record that the metric for which $J$ is orthogonal has the form 
$g=g_0(a\times a)$ on $S^6$. $\Diamond$
\vspace{0.1in}

\noindent The time has come to return to our field theory (\ref{S6YMH}). So 
consider $S^6$ with its induced orientation from the one used on ${\rm G}_2$ 
and also take the triple $(\nabla , J,g)$ constructed on $S^6$. The connection 
$\nabla$ on $TS^6$ has a corresponding curvature tensor $R_\nabla$. It is 
already meaningful to consider the symmetrized part 
${\mathfrak S}(R_\nabla )$ of this induced curvature operator on $TS^6$. 
Define the smooth function $f: S^6\rightarrow\R^+$ by $\dd V =
f\dd (y{\rm SU}(3))$ where $\dd V$ is the volume form to $g$ and 
$\dd (y{\rm SU}(3))$ denotes the standard coset measure on 
${\rm G}_2/{\rm SU}(3)\cong S^6$ induced by $\dd y=\dd V_0$ on ${\rm G}_2$. 
Inserting the equality $1=\int_{{\rm SU}(3)}\dd\gamma$ into 
(\ref{S6YMH}) and making use of the Fubini formula 
\cite[Proposition I.5.16]{bro-tom}) for the coset 
${\rm G}_2/{\rm SU}(3)\cong S^6$ we write 
\begin{eqnarray}
\ce (\nabla ,J ,g)\!\!\!&=&\!\!\!\frac{1}{2}\int\limits_{S^6}\left(
\frac{1}{e^2}\vert {\mathfrak S}(R_\nabla )\vert^2_g+
\vert\nabla g\vert^2_g+
\vert N_\nabla J\vert^2_g+e^2\vert JJ^*-
{\rm Id}_{TS^6}\vert^2_g\right)\dd V\nonumber\\
                     &=&\!\!\!\frac{1}{2}\int\limits_{S^6}\left(\:\:
\int\limits_{{\rm SU}(3)}\left(\frac{1}{e^2}\vert {\mathfrak S}(R_\nabla )
\vert^2_g+\vert\nabla g\vert^2_g+
\vert N_\nabla J\vert^2_g+e^2\vert JJ^*-
{\rm Id}_{TS^6}\vert^2_g\right)\dd\gamma\right)f\dd (y{\rm SU}(3))\nonumber\\ 
                     &=&\!\!\!\frac{1}{2}\int\limits_{{\rm G}_2}\!\!
\left(\frac{1}{e^2}\vert {\mathfrak S}(R_{\nabla_H})\vert^2_{g_H}+
\!\vert\nabla_Hg_H\vert^2_{g_H}+\!
\vert N_{\nabla_H} J_H\vert^2_{g_H}+e^2\vert
J_HJ_H^*-{\rm Id}_H\vert^2_{g_H}\right)\!\!(\pi^*f)\dd V_0\nonumber
\end{eqnarray}
consequently by the aid of (\ref{G2vakuum}) we find that
\[0\leqq\ce (\nabla ,J ,g)\leqq\Vert \pi^*f\Vert_{L^\infty ({\rm G}_2)}
\ce (\nabla_H, J_H, g_H)=0.\]
Therefore $(\nabla ,J,g)$ is a smooth vacuum solution to the 
Yang--Mills--Higgs--Nijenhuis theory (\ref{S6YMH}). 

Putting all of our findings so far together we obtain that 
$J$ is a smooth everywhere integrable almost complex structure on 
$S^6$ consequently we have arrived at the following result:

\begin{theorem}
Assume that the {\bf Principle} in Sect. \ref{two} holds. Then there exists 
a unique smooth integrable almost complex structure $J$ on $S^6$ given by the 
Higgs field in the weakly spontaneously broken vacuum $(\nabla , J, g)$ of the 
Yang--Mills--Higgs--Nijenhuis theory (\ref{S6YMH}) on $S^6$.
   
Consequently up to isomorphism there exists at least one compact complex 
manifold $X$ whose underlying real manifold is homeomorphic to the six 
dimensional sphere. $\Diamond$
\label{Xtetel}
\end{theorem}

\begin{remark}\rm 1. The new integrable almost complex structure on $S^6$ is 
nothing but the horizontal component $J_H$ of an integrable 
almost complex structure $J$ on ${\rm G}_2$ projected onto $S^6$ by a 
non-trivial twisted projection given by (\ref{B}). By \cite[p. 123]{pit} the 
moduli space of complex structures on ${\rm G}_2$ is $\C^+\sqcup\C^-$. 
However one can demonstrate\footnote{Thanks go to N.A. Daurtseva for 
calculating explicitly the whole $\C^+\sqcup\C^-$ family of integrable almost 
complex structures on ${\rm G}_2$ in the origin $e\in{\rm G}_2$ i.e., on 
${\mathfrak g}_2=T_e{\rm G}_2$.\label{natasha}} that any of these complex 
structures $J$ projects onto the same $J_H$ on the horizontal 
sub-bundle $H\subset T{\rm G}_2$ consequently the constructed structure on 
$S^6$ is unique. 

2. We have seen in Lemma \ref{homogenlemma} that although 
$Y\cong{\rm G}_2$ as a complex manifold is homogeneous since left 
translations act by biholomorphisms i.e., $[(L_y)_*,J ]=0$ for all 
$y\in {\rm G}_2$ this property does {\it not} descend to $X\cong S^6$. A 
standard way to get a homogeneous strutcture on $S^6$ is to try to push down 
$J_H: H\rightarrow H$ onto $S^6$ with respect to the canonical mapping 
$\pi_*\vert_H:H\rightarrow TS^6$. For this to happen $J_H$ should commute 
with the adjoint action of ${\rm G}_2$, see 
\cite[Volume II, Chapter X.6]{kob-nom}. However one can demonstrate by an 
explicit calculation\footnote{See Footnote \ref{natasha}.} that for $J_H$ 
this fails. In our language (cf. Definition \ref{ekv.fourier}) this 
procedure would look like this: $J_H$ on $H$ should be pushed down by the 
canonical Fourier expansion corresponsing to 
$[{\rm Id}_H]\in\cg_{{\rm SU}(3)}\backslash\cg_H$ to an endomorphism  
$J_{H/_{[{\rm Id}_H]}{\rm SU}(3)}$ of the quotient bundle  
$H/_{[{\rm Id}_H]}{\rm SU}(3)$ over ${\rm G}_2/{\rm SU}(3)\cong S^6$. Then 
the induced canonical isomorphism $H/_{[{\rm Id}_H]}{\rm SU}(3)\cong TS^6$ 
would yield a homogeneous integrable almost complex structure on $S^6$. We 
note that although this is impossible for $J_H$, it is possible to perform 
the same Fourier expansion on the other non-integrable almost complex 
structure $I_H$ on ${\rm G}_2$ with 
$[{\rm Id}_H]\in\cg_{{\rm SU}(3)}\backslash\cg_H$ . This way we just 
recover the Cayley almost complex structure $I$ on $S^6$.

Instead we are forced to push down $J_H$ from $H$ to an endomorphism 
$J_{H/_{[\alpha ]}{\rm SU}(3)}$ on the bundle $H/_{[\alpha ]}{\rm SU}(3)$ 
by a non-canonical Fourier expansion making use of a twisted action with a 
non-trivial element $[\alpha ]\in\cg_{{\rm SU}(3)}\backslash\cg_H$ as in Lemma 
\ref{vakuumlemma}. 

3. It is worth mentioning that a comparison of certain results 
offers an independent check of our Theorem \ref{Xtetel} here. 
Let $Y$ be ${\rm G}_2$ equipped with a complex structure as before. Inserting 
the decomposition $TY=V\oplus H$ into $\wedge^{p,q}(TY)$ we obtain 
\[\wedge^{p,q}(TY)\cong\wedge^{p,q}(V)\oplus\wedge^{p,q}(H)\oplus\dots\]
and let $\psi^{p,q}\in [\psi^{p,q}]\in H^{p,q}(Y;\C )$ be a 
representative of a non-zero Dolbeault cohomology class. Its 
restriction $\psi^{p,q}\vert_H$ can be cut down to $S^6$ via Fourier 
expansion as before. Let $X$ be a compact complex manifold homemorphic to 
$S^6$ with the complex structure coming from $Y$. This way one 
obtains an element $\tilde{\psi}^{p,q}\in C^\infty (\wedge^{p,q}X;\C )$. 
It may happen that $[\tilde{\psi}^{p,q}]\in H^{p,q}(X;\C )$ i.e., it 
represents a Dolbeault cohomology class.

We know the following things. On the one hand it follows 
from \cite[Proposition 4.5]{pit} that $h^{0,1}(Y)=1$ and 
$h^{0,2}(Y)=0$. On the other hand it is proved in \cite{gra2, uga} that 
$h^{0,1}(X)=h^{0,2}(X)+1$. Therefore it is suggestive to expect that 
$h^{0,1}(X)=1$ and $h^{0,2}(X)=0$. Similarly, we know from 
\cite[Proposition 4.5]{pit} that $h^{2,0}(Y)=0$, $h^{1,1}(Y)=1$, 
$h^{1,0}(Y)=0$ and $h^{1,2}(Y)=1$ meanwhile 
\cite[Proposition 3.1]{uga} 
provides us that $h^{2,0}(X)+h^{1,1}(X)=h^{1,0}(X)+h^{1,2}(X)+1$. This 
suggests that $h^{2,0}(X)=0$, $h^{1,1}(X)=1$, $h^{1,0}(X)=0$ but 
$h^{1,2}(X)=0\not=1$. However in spite of these naive considerations we 
notice that the general relationship between the Hodge numbers of $Y$ and $X$ 
is certainly not straightforward.

4. It is also worth pointing out again in the retrospective why our whole 
construction breaks down in the very similar situation of 
${\rm SO}(4n+1)/{\rm SO}(4n)\cong S^{4n}$. In our understanding the crucial 
difference between these quotients and ${\rm G}_2/{\rm SU}(3)\cong S^6$ is as 
follows. Although in the former case ${\rm SO}(4n)$ can be used to 
Fourier expand a complex structure on ${\rm SO}(4n+1)$, it fails to 
furnish $TS^{4n}$ with an almost complex structure because the 
standard representation of ${\rm SO}(4n)$ is {\it real}. Meanwhile in 
the later case ${\rm SU}(3)$ can be used not only for Fourier expansion but 
also to construct an almost complex structure on $S^6$ through its $3$ 
dimensional {\it complex} representation. This is an exceptional phenomenon 
occuring only in six dimensions and was exploited in Lemma \ref{vakuumlemma}.
\label{vegso.megjegyzesek}
\end{remark}

\noindent{\bf Acknowledgement.} The author is grateful to N.A. Daurtseva, L. 
Lempert, D. N\'ogr\'adi, E. Szab\'o, Sz. Szab\'o and R. Sz\H oke for the 
stimulating discussions. The work was partially supported by OTKA grant No. 
NK81203 (Hungary).

\end{document}